\documentclass[12pt,a4paper]{amsart}
\usepackage{amssymb,amscd}

\theoremstyle{plain}
\newtheorem{theorem}{Theorem}[subsection]
\newtheorem{lemma}[theorem]{Lemma}
\newtheorem{proposition}[theorem]{Proposition}
\newtheorem{corollary}[theorem]{Corollary}

\theoremstyle{definition}
\newtheorem{definition}[theorem]{Definition}

\newtheorem{example}[theorem]{Example}

\numberwithin{equation}{subsection}

\newcommand\bq{{\bf q}}
\newcommand\bs{{\bf s}}

\newcommand\bD{{\bf D}}
\newcommand\bG{{\bf G}}
\newcommand\bH{{\bf H}}
\newcommand\bK{{\bf K}}
\newcommand\bX{{\bf X}}

\newcommand\bbA{{\mathbb A}}

\newcommand\bbG{{\mathbb G}}
\newcommand\bbP{{\mathbb P}}

\newcommand\bbR{{\mathbb R}}
\newcommand\bbZ{{\mathbb Z}}

\newcommand\cD{\mathcal{D}}
\newcommand\cE{\mathcal{E}}
\newcommand\cF{{\mathcal F}}

\newcommand\cL{{\mathcal L}}
\newcommand\cN{{\mathcal N}}
\newcommand\cO{{\mathcal O}}
\newcommand\cS{{\mathcal S}}
\newcommand\cT{{\mathcal T}}
\newcommand\cV{{\mathcal V}}
\newcommand\cX{{\mathcal X}}

\newcommand\fa{{\mathfrak a}}

\newcommand\fg{{\mathfrak g}}
\newcommand\fh{{\mathfrak h}}

\newcommand{\hX}{\widehat{X}}

\newcommand{\tB}{\widetilde{B}}
\newcommand{\tD}{\widetilde{D}}
\newcommand{\tG}{\widetilde{G}}
\newcommand{\tH}{\widetilde{H}}
\newcommand{\tX}{\widetilde{X}}

\newcommand\tfg{\widetilde{\mathfrak g}}
\newcommand\tmu{\widetilde{\mu}}
\newcommand\trho{\widetilde{\rho}}
\newcommand\tv{\widetilde{v}}
\newcommand\tx{\widetilde{x}}
\newcommand\tvarphi{\widetilde{\varphi}}
\newcommand\tLambda{\widetilde{\Lambda}}

\newcommand\tbX{\widetilde{\bf X}}

\newcommand\tcD{\widetilde{\cD}}
\newcommand\tcV{\widetilde{\cV}}

\newcommand\ad{\operatorname{ad}}

\newcommand\diag{\operatorname{diag}}
\renewcommand\div{\operatorname{div}}

\newcommand\ord{\operatorname{ord}}

\newcommand\rk{\operatorname{rk}}

\newcommand\Aut{\operatorname{Aut}}
\newcommand\CE{\operatorname{CE}}
\newcommand\Cl{\operatorname{Cl}}
\newcommand\Cone{\operatorname{Cone}}
\newcommand\Eff{\operatorname{Eff}}
\newcommand\End{\operatorname{End}}
\newcommand\Env{\operatorname{Env}}
\newcommand\G{\operatorname{G}}
\newcommand\Hom{\operatorname{Hom}}
\newcommand\Ima{\operatorname{Im}}

\newcommand\Nef{\operatorname{Nef}}
\newcommand\Pic{\operatorname{Pic}}
\newcommand\Proj{\operatorname{Proj}}
\newcommand\PSL{\operatorname{PSL}}
\newcommand\PSO{\operatorname{PSO}}
\newcommand\SL{\operatorname{SL}}
\newcommand\SO{\operatorname{SO}}
\newcommand\Spec{\operatorname{Spec}}
\newcommand\Sp{\operatorname{Sp}}
\newcommand\Span{\operatorname{Span}}
\newcommand\Spin{\operatorname{Spin}}
\newcommand\Sym{\operatorname{Sym}}

\title{The total coordinate ring of a wonderful variety} 

\author{Michel~Brion}
\address{Universit\'e de Grenoble I\\
D\'epartement de Math\'ematiques\\
Institut Fourier, UMR 5582 du CNRS\\
38402 Saint-Martin d'H\`eres Cedex, France}
\email{Michel.Brion@ujf-grenoble.fr}

\begin{document}

\begin{abstract}
We study the cone of effective divisors and the total coordinate ring
of wonderful varieties, with applications to their automorphism group. 
We show that the total coordinate ring of any spherical variety is
obtained from that of the associated wonderful variety by a base
change of invariant subrings.
\end{abstract}

\maketitle

\tableofcontents

\section{Introduction}
\label{sec:introduction}

Wonderful varieties are smooth projective algebraic varieties where a
semi-simple algebraic group $G$ acts with finitely many orbits;
one assumes, in addition, that the $G$-orbit closures look like the
coordinate subspaces in an affine space $\bbA^r$. Examples include the
flag varieties $G/P$ (the case where $r=0$), and the complete
symmetric varieties of De Concini and Procesi (see \cite{DP83}; here
$r$ is the rank of the symmetric space embedded in the variety
as the open $G$-orbit).

The class of wonderful varieties forms a natural generalization of
that of flag varieties; it also plays a central role in the theory of 
spherical varieties, as shown by Luna. Indeed, by \cite{Lu96}, every
wonderful variety is spherical, i.e., it contains an open orbit of a
Borel subgroup of $G$. Conversely, the classification of spherical
varieties may be reduced to that of wonderful varieties, see
\cite{Lu01}. And the latter is known in many cases: when the ``rank''
$r$ is at most $2$ (see \cite{Ak83,HS82} for $r=1$, and
\cite{Wa96} for $r=2$), and when the group $G$ is of type $A$ (see
\cite{Lu01}) or $D$ (see \cite{BP05}). However, the problem of
classifying wonderful varieties is still open in general, due to the
complexity of their combinatorial invariants: the ``spherical
systems'' of \cite{Lu01}.

An important feature of this problem is that a given variety may be
wonderful with respect to several group actions, which may also have
the same orbits. For example, the projective space $\bbP^n$ is of
course a wonderful variety under $\SL_{n+1}$, and (if $n$ is odd)
under $\Sp_{n+1}$, both acting transitively. Moreover, $\bbP^n$ is a
wonderful variety under $\SO_{n+1}$ acting with two orbits: a smooth
quadric hypersurface and its complement. In addition, $\bbP^6$ is a
wonderful variety under the subgroup $\G_2 \subset \SO_7$ (acting with
the same orbits), and $\bbP^7$ is a wonderful variety under $\Spin_7$
(acting via the projectivization of its spin representation; the
orbits are still a quadric and its complement). One may show that
there is no other structure of a wonderful variety on a projective
space. Another example is $\bbP^n \times \bbP^n$ blown up along the
diagonal, which is a wonderful variety of rank $1$ under the natural
action of $\SL_{n+1}$ and also (if $n$ is odd) a wonderful variety of
rank $2$ under $\Sp_{n+1}$.

In the present paper, we study three closely related algebro-geometric
invariants of wonderful varieties: the cone of effective divisors, the
connected automorphism group, and the total coordinate ring. This
yields insight into their combinatorial invariants, and into the
classification of all wonderful group actions on a given variety.

To describe our results, recall that any wonderful variety $X$ comes
with a finite set of prime divisors called ``colors'', which yields a
basis of the Picard group. Thus, we may define (after 
\cite{BH03, Co95, EKW04}) the total coordinate ring $R(X)$ as the
direct sum of spaces of sections of all linear combinations of colors
with integral coefficients. This ring carries a natural grading by the
monoid of classes of effective divisors, $\Eff(X)$. 

It is easy to show that the algebra $R(X)$ and the monoid $\Eff(X)$
are finitely generated when $X$ is any spherical variety. For
wonderful varieties, we obtain much more precise results. Namely, we
first describe all the extremal rays of the cone $\Eff(X)_{\bbR}$ of 
effective divisors (Theorems \ref{thm:fixed} and \ref{thm:moved}). 
They fall into two types: the pull-backs of Schubert divisors under
fibrations $X \to X'$ where $X'$ is a flag variety with Picard group
$\bbZ$, and those ``boundary divisors'' (the $G$-orbit closures of
codimension $1$) that are ``fixed'', that is, any positive multiple
has a unique section up to scalars.

The fixed divisors also appear when studying automorphisms of
the wonderful $G$-variety $X$. Indeed, these divisors are all stable
under the connected automorphism group $\Aut^0(X)$. We show that 
$\Aut^0(X)$ is semi-simple, as well as any closed connected subgroup
$G'$ containing the image of $G$; moreover, $X$ is still wonderful
under the action of $G'$. In particular, $X$ has a canonical structure
of a wonderful variety under $\Aut^0(X)$; then its colors are the same
as those associated with the $G$-action, and its boundary divisors are 
exactly the fixed divisors (see Theorem \ref{thm:aut} for a more
general statement). 

In loose words, the colors and the fixed boundary divisors are
intrinsic invariants of the variety $X$. Hence the same holds for a
large part of Luna's combinatorial invariants, by their definition in
\cite{Lu01}.

We then turn to the total coordinate ring $R(X)$, also known as the 
section ring, or Cox ring, or factorial hull. The canonical sections
$s_1,\ldots,s_r$ of the boundary divisors are easily seen to form a
regular sequence in $R(X)$ and to generate the invariant subring
$R(X)^G$; moreover, all fibers of the corresponding flat morphism 
$q : \Spec R(X) \to \bbA^r$ are normal varieties (Proposition
\ref{prop:quot}).

We also obtain an algebraic descriptions of $R(X)$ as a certain
``Rees algebra'' (Theorem \ref{thm:rees}), and a presentation by
generators and relations. It should be emphasized that the simple
description of the total coordinate ring of a toric variety (a
polynomial ring over the coordinate ring of the acting torus, see
\cite{Co95, AHS02}) does not extend in a straightforward way
to a wonderful variety $X$. Specifically, the ring $R(X)$ still admits
simple generators, namely, $s_1,\ldots,s_r$ and the $G$-submodules
spanned by the canonical sections of colors, but their relations turn
out to be rather complicated.

For example, $R(G/P)$ is the coordinate ring of the quasi-affine
homogeneous space $G/[P,P]$, which is defined by quadratic
equations. In the general case, the quotient $R(X)/(s_1,\ldots,s_r)$
is also defined by quadratic equations, that one has to lift to $R(X)$
(Proposition \ref{prop:rel}). This yields some information on the
(still largely open) problem of describing the product of any two
simple $G$-submodules of $R(X)$; see Proposition \ref{prop:gs}. For
certain complete symmetric varieties, more precise results are due to 
Chiriv\`{\i}, Littelmann and Maffei, see \cite{CLM06}.

Finally, we study the total coordinate ring of an arbitrary spherical
variety $X$ under a connected reductive group $G$. Here $X$ may be
quite singular, so that it is natural to define the ring $R(X)$ via
Weil divisors as in \cite{BH03}; actually, one should consider
$G$-linearized divisorial sheaves to keep track of the $G$-action. 
Then the Picard group is replaced with the equivariant divisor class
group $\Cl^G(X)$, which has no canonical basis in general; in fact,
this group may have torsion. So the definition of $R(X)$ requires
some work, see \ref{subsec:total} for details.

Once this is done, the ring $R(X)$ turns out to be obtained from the
total coordinate ring of a wonderful variety by a base
change. Specifically, the spherical $G$-variety $X$ is equipped with 
an equivariant rational map $\varphi: X - \to \bX$, where $\bX$ is
a wonderful variety under a semi-simple quotient $\bG$ of the group
$G$; see \cite{Lu01}. This yields an equivariant homomorphism 
$\varphi^*: R(\bX) \to R(X)$; we show that the induced map 
$R(\bX)\otimes_{R(\bX)^{\bG}} R(X)^{\bG} \to R(X)$ is an isomorphism 
(Theorem \ref{thm:cart}). This may be regarded as a universal property
of the morphism $q : \Spec R(\bX)^{\bG} \to \bbA^r$, a flat family of
affine $\bG$-varieties. The corresponding morphism to the 
``invariant Hilbert scheme'' of \cite{AB05} calls for investigation. 

We also study the connected automorphism group $\Aut^0(X)$, where $X$
is an arbitrary complete toroidal variety (i.e., the above rational
map $\varphi$ is a morphism). In contrast with the wonderful case,
$\Aut^0(X)$ is generally far from being semi-simple (e.g., it is an
$n$-dimensional torus when $X$ is the projective $n$-space blown up at
$n+1$ general points) or even reductive (e.g., for Hirzebruch
surfaces). However, we show that the ``double centralizer'' of $G$
(that is, the centralizer in $\Aut^0(X)$ of the equivariant connected
automorphism group $\Aut^0_G(X)$) is reductive, as well as any 
subgroup containing the image of $G$. Also, the subgroup of
$\Aut^0(X)$ preserving the boundary is reductive and depends only on
the open $G$-orbit in $X$. For this, we relate these subgroups to the
corresponding subgroups of $\Aut^0(\bX)$, see Theorem \ref{thm:aut-bis}.

For a complete symmetric variety $X$, the ring $R(X)$ has been
considered by Chiriv\`{\i} and Maffei \cite{CM03} who constructed a 
``standard monomial theory'' in this setting. This ring also appears,
somewhat implicitly, in Vinberg's work \cite{Vi95b} on reductive
algebraic semigroups which has been the starting point of the present
article. Indeed, as shown by Rittatore in \cite{Ri97,Ri02}, Vinberg's 
``enveloping semigroup'' $\Env(G)$ is the affine variety with
coordinate ring $R(\bar{G}_{\ad})$, where $\bar{G}_{\ad}$ denotes 
the wonderful compactification of the adjoint group $G_{\ad}$ (a
special case of a wonderful symmetric variety); see \cite{Re06} 
and its references for related results. Throughout this paper,
we shall specialize our results to the group compactifications
$\bar{G}_{\ad}$ in a series of examples.

\medskip

\noindent
{\it Acknowledgments.} A large part of this work was written in
January 2006, during a staying at the University of Georgia. I would
like to thank this institution for its support; also, many thanks are
due to Valery Alexeev and Bill Graham for hospitality and fruitful
discussions. Finally, I thank the referee for his careful reading and
valuable comments.

\medskip

\noindent
{\bf Notation.}
We shall use \cite{Ha77} as a general reference for algebraic
geometry, and \cite{PV94} for algebraic transformation groups.

The ground field $k$ is algebraically closed, of characteristic zero. 
By a variety, we mean an integral separated scheme of finite type over
$k$. We denote by $k(X)$ the field of rational functions on a
variety $X$, and by $k[X] = \Gamma(X,\cO_X)$ its algebra of global
regular functions.

We denote by $G$ a connected reductive linear algebraic group; we
choose a Borel subgroup $B\subset G$ and a maximal torus 
$T\subset B$. Let $W := N_G(T)/T$ be the Weyl group, and $U$ the
unipotent part of $B$.
The group of multiplicative characters $\cX(T)$ is also denoted by
$\Lambda$, and identified with the character group $\cX(B)$. Let 
$\Phi \subset \Lambda$ be the root system of $(G,T)$ with its subset
of positive roots $\Phi^+$ consisting of roots of $(B,T)$. This
defines the subset $\Lambda^+ \subset \Lambda$ of dominant weights,
which is also the set of highest weights of simple rational
$G$-modules. For any $\lambda \in \Lambda$, we denote by $V(\lambda)$
the corresponding simple module.

\section{Divisors on wonderful varieties}
\label{sec:divisors}

In this section, $G$ is assumed to be semi-simple and simply-connected.

\subsection{Wonderful varieties}\label{subsec:wonderful} 

Recall from \cite{Lu96} that a $G$-variety $X$ is called 
{\it wonderful of rank $r$} if it satisfies the following conditions: 

\smallskip

\noindent
(1) $X$ is complete and non-singular.
\smallskip

\noindent
(2) $X$ contains an open $G$-orbit $X^0$, and the complement 
$X \setminus X^0$ (the \emph{boundary}) is a union of $r$ non-singular
prime divisors $X_1,\ldots,X_r$ (the {\it boundary divisors}) with
normal crossings.

\smallskip

\noindent
(3) The $G$-orbit closures in $X$ are exactly the partial intersections
$\bigcap_{i \in I} X_i$, where $I$ is a subset of $\{1,\ldots,r\}$.

\smallskip

In particular, $X$ contains a unique closed $G$-orbit, the
intersection of all boundary divisors; it follows that $X$ is
projective. We now recall a number of notions and results concerning
wonderful varieties; see \cite{Lu01} and its references for details
(some of our notation differs from \cite{Lu01}.)

Denote by $X^0_B$ the open $B$-orbit in the wonderful variety
$X$. Then $X^0_B$ is affine and contained in $X^0$, so that 
$X^0 \setminus X^0_B$ is a union of finitely many prime $B$-stable
divisors, called {\it colors}. Let $\cD$ be the set of colors,
regarded as divisors in $X$ by taking closures; these are the prime
$B$-stable divisors in $X$ which are not $G$-stable. Moreover,
\begin{equation}
X \setminus X^0_B = X_1 \cup \cdots \cup X_r \cup
\bigcup_{D \in \cD} D.
\end{equation}
Let also 
$$
X_B := X \setminus \bigcup_{D \in \cD} D,
$$
this is an open $B$-stable subset of $X$, called the 
\emph{big cell}. In fact, $X_B$ is isomorphic to an affine space, and
meets each $G$-orbit along its open $B$-orbit. 

Each intersection $X_i\cap X_B$ has an equation $f_i \in k[X_B]$,
uniquely defined up to a non-zero scalar. In particular, $f_i$ is a
$B$-eigenvector; let $\gamma_i \in \Lambda$ be the opposite of its
weight. Then $\gamma_1, \ldots, \gamma_r$ are the 
{\it spherical roots} of $X$. In fact, the spherical roots form 
a basis of a reduced root system $\Phi_X$, and the corresponding
Weyl group $W_X$ is a subgroup of $W$. Moreover, each spherical root
is either a positive root of $\Phi$, or a sum of two such roots.

We may regard $f_1,\ldots,f_r$ as $B$-eigenvectors in $k(X)$; in fact,
the multiplicative group of $B$-eigenvectors $k(X)^{(B)}$ is freely
generated by $f_1,\ldots,f_r$ over its subgroup $k^*$ of
constants. Moreover, associating with each eigenvector its weight
yields an injective homomorphism $k(X)^{(B)}/k^* \to \Lambda$; its
image is the {\it weight group} $\Lambda(X)$, freely generated by the
spherical roots.

Choose a base point $x \in X^0_B$ and denote by $H = G_x$ its isotropy
group. This yields isomorphisms $X^0 \cong G/H$, 
$X^0_B \cong BH/H \cong B/B\cap H$, and  
$\Lambda(X) \cong \cX(B)^{B\cap H}$ (the group of characters of $B$
which are invariant under $B\cap H$). We may choose $T$ so that
$T \cap H$ is a maximal diagonalizable subgroup of $B \cap H$; then 
$\cX(B)^{B\cap H} \cong \cX(T/T\cap H)$. 

Also, recall that $H$ has finite index in its normalizer $N_G(H)$;
moreover, the equivariant automorphism group $\Aut_G(X)$ is isomorphic
to $\Aut_G(G/H) \cong N_G(H)/H$ via restriction. In particular, 
$\Aut_G(X)$ is finite.

Let $\pi:G \to G/H$ denote the projection. Each $D \in \cD$ yields a
$B \times H$-stable divisor $\pi^{-1}(D)$ on $G$, which has an equation 
$f_D \in k[G]$ uniquely defined up to a non-zero scalar (since $k[G]$
is a UFD and $k[G]^* = k^*$). Thus, $f_D$ is an eigenvector of 
$B \times H$, and $f_D(1) \neq 0$ so that $f_D$ is uniquely determined
by requiring that $f_D(1) = 1$. Likewise, each $f_i$ is uniquely
determined by requiring that $f_i(x)=1$. 

By \cite[Lem.~6.2.2]{Lu01}, the $B\times H$-eigenvectors in $k(G)$ are
exactly the Laurent monomials in the $f_D$, $D \in \cD$; it follows
that the $B\times H$-eigenvectors in $k[G]$ are the monomials in the
$f_D$'s. Let $(\omega_D,\chi_D)  \in \cX(B) \times \cX(H)$ be the
weight of $f_D$. Note that the restrictions to $B\cap H$ of $\omega_D$
and $\chi_D$ coincide. In other words, $(\omega_D,\chi_D)$ is in the 
fibered product $\cX(B) \times_{\cX(B\cap H)} \cX(H)$. In fact, 
we have the following useful result (probably known, but for which 
we could not supply a reference): 

\begin{lemma}\label{lem:str}
The abelian group $\cX(B) \times_{\cX(B\cap H)} \cX(H)$
is freely generated by the pairs $(\omega_D,\chi_D)$,
$D \in \cD$. Moreover, $\cX(H)$ is generated by the $\chi_D$, 
$D \in \cD$.
\end{lemma}

\begin{proof}
Consider the action of the group $B \times H$ on $G$ via 
$(b,h) \cdot g = bgh^{-1}$. The product $BH$ is the orbit of the
identity element $1$; it is open in $G$ (since $B/B\cap H$ is open in
$G/H$), and the isotropy group of $1$ equals $\diag(B\cap H)$. It
follows that each $B\times H$-eigenvector in 
$k(G) = k(BH) = k((B\times H)/\diag(B\cap H))$ 
is determined up to a scalar by its pair of weights. Moreover,
the set of these weights is $\cX(B) \times_{\cX(B\cap H)} \cX(H)$. 
This implies the first assertion. For the second one, observe that the 
projection $\cX(B) \times_{\cX(B\cap H)} \cX(H) \to \cX(H)$
is surjective, since any character of $B\cap H$ extends to a character
of $B$.
\end{proof}

Any prime divisor $D$ in $X$ defines a discrete valuation $v_D$ of the
function field $k(X)$; its restriction to $k(X)^{(B)}$ yields a
homomorphism 
$$
\rho(v_D): k(X)^{(B)}/k^* \cong \Lambda(X) \to \bbZ~.
$$
Let $v_i := \rho(v_{X_i})$ for $i=1,\ldots,r$, then 
$v_i(f_j) = \delta_{ij}$. In other words, $v_1,\ldots,v_r$ 
form a basis of the abelian group $\Hom(\Lambda(X),\bbZ)$, dual to the 
basis $(-\gamma_1,\ldots,-\gamma_r)$. In particular, 
$\rho(v_D) = \sum_{i=1}^r \langle D,\gamma_i \rangle \, v_i,$
where 
\begin{equation}\label{eqn:def}
\langle D,\gamma_i\rangle := - \rho(v_D)(f_i).
\end{equation}
Also, $f_i = \prod_{D \in \cD} f_D^{- \langle D,\gamma_i \rangle}$. 
Equivalently,
\begin{equation}\label{eqn:weights}
\sum_{D \in \cD} \langle D,\gamma_i \rangle \, \omega_D = \gamma_i
\quad \text{and} \quad
\sum_{D \in \cD} \langle D,\gamma_i \rangle \, \chi_D =0.
\end{equation}
The convex cone
$$
\cV := \Cone(v_1,\ldots,v_r) \subset \Hom(\Lambda(X),\bbR)
$$
may be identified with the set of $G$-invariant discrete valuations of
$k(X) = k(G/H)$. Thus, $\cV$ is called the {\it valuation cone} of the
spherical homogeneous space $G/H$, and $\Lambda(X) = \Lambda(G/H)$ is
called its \emph{weight group}. 

The cone $\cV$ and the map $\rho : \cD \to \Hom(\Lambda(G/H),\bbR)$ are
the main ingredients of the Luna--Vust classification of embeddings of
$G/H$ and, more generally, of any spherical homogeneous space. We
refer to \cite{Kn91} for an exposition of this theory; we shall freely
use some of its results in the sequel. We shall also need the
following 

\begin{lemma}\label{lem:cones}
The opposite of the valuation cone is contained in the cone generated
by the $\rho(v_D)$, $D \in \cD$.
\end{lemma}
  
\begin{proof}
It suffices to show the inclusion of dual cones in $\Lambda(X)_{\bbR}$:
$$
\Cone(\rho(v_D), D \in \cD)^{\vee} \subseteq -\cV^{\vee} 
= \Cone(\gamma_1,\ldots,\gamma_r).
$$
Let $f \in k(X)^{(B)}$ such that $\rho(v_D)(f) \geq 0$ for all 
$D \in \cD$. Then $f$ is a regular function on $G/H$, and hence its
weight $\lambda$ is dominant. Thus,
$(\lambda,\gamma_i) \geq 0$ for $i = 1,\ldots,r$, where $(\,,)$ denotes
the scalar product on $\Lambda(X)_{\bbR} \subseteq \Lambda_{\bbR}$
defined by the Killing form of $G$. Since $(\,,)$ is invariant under the
Weyl group $W_X$ of the root system $\Phi_X$, it follows that
$\lambda$ is also dominant relatively to this root system and its
basis of spherical roots. Thus, the coordinates of $\lambda$
in this basis are non-negative, as follows from the equality
$\vert W_X \vert \, \lambda = \sum_{w \in W_X} (\lambda - w(\lambda))$.
\end{proof}
 
\begin{example}\label{ex:group1}
Let $Z$ denote the center of $G$, and $G_{\ad} := G/Z$
the corresponding adjoint group. Then $G_{\ad}$ is a homogeneous
space under $G\times G$ acting by left and right multiplication; the
isotropy group of the base point $1$ is $(Z \times Z)\diag(G)$. 
By \cite{DP83}, $G_{\ad}$ admits a wonderful compactification
$X = \bar{G}_{\ad}$, of rank $r = \rk(G)$.  Its spherical roots
(relatively to the Borel subgroup $B^- \times B$, where $B^-$ denotes
the unique Borel subgroup such that $B^- \cap B = T$) are the pairs
$(-\alpha_i,\alpha_i)$ where $\alpha_1,\ldots,\alpha_r$ are the simple
roots of $\Phi$. Thus, $\Lambda(X)$ may be identified with the root
lattice of $\Phi$; then $\Hom(\Lambda(X),\bbZ)$ is the weight lattice,
and $v_1,\ldots,v_r$ are the opposites of the fundamental
co-weights. The valuation cone is the negative Weyl chamber. The map   
$\rho : \cD \to \Hom(\Lambda(X),\bbZ)$ is injective, and its image
consists of the simple co-roots
$\alpha_1^{\vee},\ldots,\alpha_r^{\vee}$. Hence the integers 
$\langle D,\gamma_i \rangle$ are the coefficients of the Cartan matrix
of the dual root system $\Phi^{\vee}$. 
\end{example}

\subsection{Picard group}
\label{subsec:picard}

Since the wonderful variety $X$ is non-singular, we may identify the
Picard group $\Pic(X)$ with the divisor class group. And since the
complement of the union of the colors is an affine space, the classes
of these colors form a basis of $\Pic(X)$. We shall identify $\Pic(X)$
with the free abelian group $\bbZ^{\cD}$ and denote its elements by
$\sum_{D \in \cD} n_D [D]$, where $(n_D) \in \bbZ^{\cD}$. 

Recall from \cite[Sec.~2.6]{Br89b} that the divisor
$\sum_{D \in \cD} n_D D$ is base-point-free (resp.~ample) if and
only if $n_D \geq 0$ (resp.~$>0$) for all $D \in \cD$. In other words,
the monoid $\Nef(X)$ consisting of the classes of numerically
effective (nef) divisors is generated by the classes of colors, and
each nef divisor is base-point-free. In particular, the nef cone
$\Nef(X)_{\bbR}$ is simplicial.

Since $G$ is semi-simple and simply-connected, each invertible sheaf
on $X$ admits a unique $G$-linearization, see \cite[Sec.~2.4]{KKLV}. Thus,
we shall also identify $\Pic(X)$ with the group of isomorphism classes
of $G$-linearized invertible sheaves, $\Pic^G(X)$. We shall use the
notation $\cO_X(\sum_{D \in \cD} n_D D)$
when considering an element of $\Pic(X)$ as a $G$-linearized
invertible sheaf. Each space of global sections
$\Gamma(X,\cO_X(\sum_{D \in \cD} n_D D))$ has a canonical $G$-module
structure.

In particular, each effective divisor admits a $B$-semi-invariant
global section. In other words, the monoid $\Eff(X)$ consisting of the
classes of effective divisors is generated by the classes of colors
and of boundary divisors. To express the latter in terms of the
former, note that
\begin{equation}\label{eqn:idiv}
\div(f_i) =  X_i - \sum_{D \in \cD} \langle D,\gamma_i \rangle \, D
\end{equation}
by (\ref{eqn:def}), so that 
\begin{equation}\label{eqn:irel}
[X_i] = \sum_{D \in \cD} \langle D,\gamma_i \rangle \, [D]
\quad \text{in } \Pic(X).
\end{equation}

We now describe the restriction map 
$i^*: \Pic(X) \to \Pic(G/H)$, where $i: G/H \to X$ denotes the
inclusion. Recall the isomorphisms
$$
\Pic(G/H) \cong \Pic^G(G/H) \cong \cX(H),
$$
where the latter is obtained by restriction to the base point. The
inverse isomorphism is given by 
\begin{equation}\label{eqn:ass}
\cX(H) \to \Pic^G(G/H), \quad \chi \mapsto \cL(\chi),
\end{equation}
where $\cL(\chi)$ denotes the $G$-linearized invertible sheaf on $G/H$ 
associated with the character $\chi$. Moreover, the $G$-module 
$\Gamma(G/H, \cL(\chi))$ is isomorphic to $k[G]^{(H)}_{\chi}$
(the space of right $H$-eigenvectors of weight $\chi$); see, e.g.,
\cite[Sec.~2.1]{KKLV} for these facts. 

\begin{proposition}\label{prop:sequ}
(1) The sequence
\begin{equation}\label{eqn:res}
\CD
0 @>>>  \bbZ^r @>{u}>>  \Pic(X) @>{i^*}>> \Pic(G/H) @>>> 0
\endCD
\end{equation}
is exact, where $u(x_1,\ldots,x_r) := \sum_{i=1}^r x_i \, [X_i]$.

\smallskip

\noindent
(2) The canonical section $s_D \in \Gamma(X,\cO_X(D))$ is a
$B$-eigenvector of weight $\omega_D$; the canonical section 
$s_i \in \Gamma(X,\cO_X(X_i))$ is $G$-invariant.

\smallskip

\noindent
(3) The natural map $\cX(B) \to \Pic^B(X_B)$ is an isomorphism. Via
the inverse isomorphism, the restriction of $\cO_X(D)$
(resp. $\cO_X(X_i)$) is mapped to $\omega_D$ (resp. $\gamma_i$).

\smallskip

\noindent
(4) The product of restriction maps
$$
\pi: \Pic(X) = \Pic^G(X) \to \Pic^B(X_B) \times \Pic^G(G/H) 
\cong \cX(B) \times \cX(H)
$$ 
is injective, and its image is 
$\cX(B) \times_{\cX(B\cap H)} \cX(H)$.
Together with the isomorphism (\ref{eqn:ass}), this
identifies (\ref{eqn:res}) with the exact sequence
\begin{equation}\label{eqn:char}
0 \to \Lambda(X) = \cX(B)^{B\cap H} = \bigoplus_{i=1}^r \bbZ \gamma_i 
\to \bbZ^{\cD} = \cX(B) \times_{\cX(B\cap H)} \cX(H) 
\to \cX(H) \to 0,
\end{equation}
where the map on the right is the second projection. Moreover, each
$[X_i]$ is identified with the spherical root $\gamma_i$.
\end{proposition}

\begin{proof}
(1) Since $G/H = X \setminus \bigcup_{i=1}^r X_i$, the map $i^*$ is
surjective and its kernel is generated by the classes
$[X_1], \ldots, [X_r]$. They are linearly independent: if $f \in k(X)$
has zeroes and poles along $X_1,\ldots,X_r$ only, then $f$ must be
$G$-semi-invariant, and hence constant.

(2) The restriction of $s_D$ to $G/H$ equals 
$$
f_D \in k[G]^{(H)}_{\chi_D} \cong \Gamma(G/H,\cL(\chi_D)) 
= \Gamma(G/H,\cO_X(D)),
$$
a $B$-eigenvector of weight $\omega_D$. The $G$-invariance of $s_i$
is readily checked.

(3) The first assertion follows from the fact that $k[X_B]$ is a UFD
and $k[X_B]^* = k^*$. For the second assertion, observe that $s_D$
yields a trivialization of $\cO_X(D)$ over $X_B$. Likewise, $f_i$
yields a trivialization of $\cO_X(X_i)$ over $X_B$.
 
(4) $\pi$ maps any class $[D]$ to the pair 
$(\omega_D,\chi_D)$, and any class $[X_i]$, to
$(\gamma_i,0)$. Together with Lemma \ref{lem:str}, this implies the  
assertions.
\end{proof}

Next we determine the sections of the positive multiples of boundary
divisors.

\begin{definition}
We say that a prime divisor $D$ in $X$ is \emph{fixed} if the  
global sections of any positive integral multiple $nD$ are just the
scalar multiples of the canonical section $s_{nD} = s_D^n$. 
Equivalently, any regular function on $X \setminus D$ is constant. 
\end{definition}

\begin{theorem}\label{thm:fixed}
Each boundary divisor $X_i$ is either fixed or base-point-free.
In the latter case, the spherical root $\gamma_i$ is a dominant
weight, and there is an isomorphism of $G$-modules
\begin{equation}\label{eqn:dec}
\Gamma(X,\cO_X(nX_i)) \cong \bigoplus_{m=0}^n V(m\gamma_i).
\end{equation}
Moreover, the $G$-variety
$$
X' := \Proj \bigoplus_{n \geq 0} \Gamma(X,\cO_X(nX_i))
$$
is a wonderful variety of rank $1$; its boundary divisor $X'_i$ is
ample, and $X_i = \varphi^{-1}(X'_i)$ where $\varphi : X \to X'$ is
the natural morphism.
\end{theorem}

\begin{proof}
Let $f \in k(X)^{(B)}$ such that the divisor $\div(f) + n X_i$ is
effective. Equivalently, $v_D(f) \geq 0$ for all $D \in \cD$, 
$v_j(f) \ge 0$ for all $j \neq i$, and $v_i(f) + n \geq 0$. By Lemma
\ref{lem:cones}, it follows that $v_j(f) = 0$ for all $j \ne i$, that
is, $f$ is a scalar multiple of some power $f_i^{-m}$, where 
$m \in \bbZ$. Then, by (\ref{eqn:idiv}),
$$
\div(f) + n X_i = (n-m) X_i + m 
\sum_{D \in \cD}  \langle D,\gamma_i \rangle \, D.
$$
Moreover, there exists $D \in \cD$ such that 
$\langle D,\gamma_i \rangle > 0$, since $X_i$ cannot be
linearly equivalent to the opposite of an effective divisor. 
Thus, $0 \leq m \leq n$.

If there exists $D \in \cD$ such that $\langle D,\gamma_i \rangle <0$,
then $m=0$, i.e., $f$ is constant: $X_i$ is fixed.

Otherwise, $X_i$ is base-point-free, 
$\gamma_i = \sum_{D \in \cD} \langle D,\gamma_i \rangle \, \omega_D$
is dominant, and the $B$-eigenvectors in $\Gamma(X,\cO_X(nX_i))$ are
exactly the scalar multiples of the $s_i^n f_i^{-m}$, where 
$0 \leq m \leq n$. This implies (\ref{eqn:dec}). As a consequence, the
$B$-eigenvectors in the $G$-algebra 
$\bigoplus_{n \geq 0} \Gamma(X,\cO_X(nX_i))$ are exactly the monomials
in $s_i$ and $s_i f_i^{-1}$.

By construction, $X'$ is a normal projective $G$-variety equipped with
a $G$-equivariant morphism $\varphi: X \to X'$ such that the induced map 
$\cO_{X'} \to \varphi_* \cO_X$ is an isomorphism; in particular,
$\varphi$ is surjective with connected fibers. Moreover, 
$X_i = \varphi^{-1}(X'_i)$ where $X'_i$ is an ample $G$-stable divisor
on $X'$. Together with the preceding argument, it follows that $X'$ is
spherical of rank $1$ with boundary $X'_i$; hence $X'$ is wonderful.
\end{proof}

The $G$-varieties $X'$ that appear in Theorem \ref{thm:fixed}
are precisely the wonderful $G$-varieties of rank $1$ having an ample 
boundary divisor. Equivalently, their open orbit $G/H'$ is affine
(indeed, the complement of an ample divisor is affine; conversely, if
$G/H'$ is affine, then its complement is ample as follows, e.g., from
Theorem \ref{thm:fixed}); in other words, the isotropy group $H'$ is
reductive. These varieties have been classified by Akhiezer
\cite[Thm.~4]{Ak83} over the complex numbers; he showed that they
correspond bijectively to compact isotropic Riemannian
manifolds (see also \cite{Br89a} for a more algebraic approach). It
turns out that the varieties $X'$ are all  homogeneous under their
full automorphism group; this is deduced in \cite{Ak83} from the
classification, while an a priori proof is obtained in
\cite[Cor.~4.1.3]{BB96}. We shall show in Theorem \ref{thm:moved} that
$\varphi: X \to X'$ is a homogeneous fibration with respect to an
overgroup of $G$. 

\begin{example}\label{ex:group2} 
Let $X = \bar{G}_{\ad}$ as in Example \ref{ex:group1}. Then
the Picard group may be identified with the weight lattice of
$G$. The classes of colors (resp.~boundary divisors) are then
identified with the fundamental weights (resp.~simple roots). The nef
cone is the positive Weyl chamber. The fixed boundary divisors
correspond to those simple roots $\alpha_i$ that are dominant.
Equivalently, $\alpha_i$ is orthogonal to all other simple roots,
that is, the image of the associated morphism $\SL_2 \to G$ yields a
direct factor of $G_{\ad}$ isomorphic to $\PSL_2$. Then the variety
$X'$ of Theorem \ref{thm:fixed} is the wonderful compactification
of $\PSL_2$, i.e., the projectivization $\bbP(M_2)$ of the space of 
$2 \times 2$ matrices; it is a direct factor of $X$, and the map
$\varphi$ is just the corresponding projection.
\end{example}

\subsection{The cone of effective divisors}
\label{subsec:cone}

Denote by $\Eff(X)_{\bbR}$ the cone in $\Pic(X)_{\bbR} \cong \bbR^{\cD}$
generated by the classes of effective divisors. This 
\emph{cone of effective divisors} is generated by the classes of the
colors and of the boundary divisors; in particular, it is a rational
polyhedral convex cone, which contains no line as $X$ is complete.  
The inclusion $\Eff(X) \subseteq \Eff_{\bbR}(X) \cap \Pic(X)$ may be
strict; in other words, a non-effective divisor may have an
effective non-zero multiple, see Example \ref{ex:group3} below.

We first describe the intersection of the cone of effective
divisors with the subspace spanned by the classes of boundary
divisors (these are linearly independent by Proposition
\ref{prop:sequ}).

\begin{lemma}\label{lem:inv}
$\Eff(X)_{\bbR} \cap \Span([X_1], \ldots, [X_r])
= \Cone([X_1],\ldots,[X_r])$.
\end{lemma}

\begin{proof}
It suffices to show that if $x_1,\ldots,x_r$ are integers such that
$\sum_{i=1}^r x_i \, [X_i]$ is effective, then they are all nonnegative.
There exists $f \in k(X)^{(B)}$ such that
$$
\sum_{i=1}^r x_i \, X_i = \div(f) + \sum_{i=1}^r y_i \, X_i 
+ \sum_{D \in \cD} z_D \, D,
$$
where $y_i,z_D \geq 0$. So $v_D(f) \leq 0$ for all $D$, and hence
$v_i(f) \geq 0$ for all $i$ by Lemma \ref{lem:cones}. Hence
$x_i = v_i(f) + y_i \geq 0$.
\end{proof}

We now turn to the description of the extremal rays of
$\Eff(X)_{\bbR}$. Clearly, any such ray is generated by the class
of a boundary divisor or a color; moreover, the classes of colors are
pairwise non-proportional, and the same holds for the classes of
boundary divisors. But a boundary divisor may well be proportional to
a color, e.g., when $\Pic(X)$ has rank $1$. So we may obtain three
types of extremal rays, considered separately in the following lemmas.

\begin{lemma}\label{lem:boundary}
For a boundary divisor $X_i$, the following conditions are equivalent: 

\begin{enumerate}

\item
$[X_i]$ generates an extremal ray of $\Eff(X)_{\bbR}$ which does
not contain the class of a color. 

\item
$X_i$ is fixed.

\item
$\langle D,\gamma_i \rangle <0$ for some $D \in \cD$.

\end{enumerate}

\end{lemma}

\begin{proof}
(1) is equivalent to the non-existence of  
$f \in k(X)^{(B)}\setminus k^*$ such that $\div(f) + n X_i $ is 
effective for some positive integer $n$. In turn, this is equivalent to 
$k[X \setminus X_i]^{(B)} = k$. But since $k[X \setminus X_i]$ is a
rational $G$-module, the latter condition amounts to 
$k[X \setminus X_i] = k$. Thus, (1) $\Leftrightarrow$ (2).

(2) $\Leftrightarrow$ (3) follows from (\ref{eqn:irel}) and
Theorem \ref{thm:fixed}.

\end{proof}

\begin{lemma} \label{lem:color}
For a color $D$, the following conditions are equivalent:

\begin{enumerate}

\item
$[D]$ generates an extremal ray of $\Eff(X)_{\bbR}$ which does not
contain the class of a boundary divisor.

\item
There exists a $G$-equivariant morphism 
$$
\varphi : X \to G/P
$$ 
where $P\supseteq H$ is a maximal parabolic subgroup of $G$, such that
$D$ is the preimage of the Schubert divisor (the unique $B$-stable
prime divisor) in $G/P$.

\end{enumerate}

\end{lemma}

\begin{proof}
As in the proof of Lemma \ref{lem:boundary}, (1) is equivalent to the
condition that any $B$-eigenvector in $\Gamma(X,\cO_X(nD))$ is a scalar
multiple of $s_D^n$, for all $n\ge 1$. In turn, this is equivalent to  
the simplicity of each $G$-module $\Gamma(X,\cO_X(nD))$, and hence to
the condition that the $G$-variety
$$
X' := \Proj \bigoplus_{n \ge 0} \Gamma(X,\cO_X(nD))
$$
is isomorphic to some $G/P$. Since $D$ is base-point-free, $X$ comes
with a $G$-morphism $\varphi : X \to X'$, and hence we may assume that
$P \supseteq H$. Moreover, $\varphi(D)$ is a prime $B$-stable ample
divisor in $G/P$, so that the parabolic subgroup $P$ is maximal.
\end{proof}

\begin{lemma}\label{lem:mixed}
For a boundary divisor $X_i$ and a color $D$, the following conditions
are equivalent:

\begin{enumerate}

\item
$[X_i]$ is a multiple of $[D]$.

\item
There exists a surjective $G$-equivariant morphism
$$
\varphi : X \to X'
$$ 
with connected fibers, where $X'$ is a wonderful $G$-variety of rank
$1$ having a unique color $D'$. Moreover, $D = \varphi^{-1}(D')$ and
$X_i = \varphi^{-1}(X'_i)$, where $X'_i$ denotes the unique boundary
divisor of $X'$.

\end{enumerate}

If (1) or (2) holds, then $[X_i]$ generates an extremal ray of
$\Eff(X)_{\bbR}$. 

\end{lemma}

\begin{proof}
(1) $\Leftrightarrow$ (2) follows from Theorem \ref{thm:fixed}.
It also follows that $\div(f_i) = X_i - aD$ for some positive integer
$a$, and that the effective $B$-stable divisors linearly equivalent to
a given positive multiple $nX_i$ are exactly the 
$(n-m)X_i + maD$, where $m=0,\ldots,n$. This implies the final assertion.
\end{proof}

The varieties $X'$ that appear in Lemma \ref{lem:mixed} are exactly
the wonderful varieties of rank $1$ and Picard group $\bbZ$ (in
particular, their boundary divisor is ample). Their classification
follows from that in \cite[Thm.~4]{Ak83}.

\begin{example}\label{ex:group3}
If $X = \bar{G}_{ad}$, then the extremal rays of $\Eff(X)_{\bbR}$
are generated by the simple roots $\alpha_1,\ldots,\alpha_r$, while
those of $\Nef(X)_{\bbR}$ are generated by the fundamental weights 
$\omega_1,\ldots,\omega_r$. 

If $G$ is simple of type $B_2$, then 
$\omega_1 - \omega_2 = \frac{1}{2}\alpha_1$, so that the colors
$D_1,D_2$ satisfy $[D_1] - [D_2] = \frac{1}{2} [X_1]$.
Thus, $E := D_1 - D_2$ is a non-effective (integral) divisor such 
that $2E$ is effective.
\end{example}

Returning to an arbitrary wonderful $G$-variety $X$, if the isotropy
group $H$ is not contained in any proper parabolic subgroup of $G$,
then the cone $\Eff(X)_{\bbR}$ is generated by the classes of boundary
divisors by Lemma \ref{lem:color}. Together with Proposition
\ref{prop:sequ}, it follows that $\Eff(X)_{\bbR}$ is simplicial. But
this need not hold in general, as shown by:

\begin{example}
In the product 
$$
\bbP^2 \times \bbP^2 \times \check \bbP^2 \times \check \bbP^2
$$
where $\check \bbP^2$ denotes the projective plane of lines in $\bbP^2$, 
let $X$ be the closed subvariety defined by the incidence conditions 
$$
p_1 \in d_1, \; p_2 \in d_1, \; p_2 \in d_2
$$
where $p_1,p_2 \in \bbP^2$ and $d_1,d_2 \in \check \bbP^2$.
Then $X$ is stable under the diagonal action of $G := \SL_3$ and one
easily checks that $X$ is a wonderful $G$-variety with boundary
divisors $X_1$ (the locus where $p_1 = p_2$) and $X_2$ 
(where $d_1 = d_2$).
Moreover, $X$ has $4$ colors, the pull-backs of the $B$-stable
lines in $\bbP^2$, $\check \bbP^2$ via the $4$ projections. 
One also checks (either directly or by using the above results) that
the cone $\Eff(X)_{\bbR}$ in $\Pic(X)_{\bbR} \cong \bbR^4$ has $6$
extremal rays, corresponding bijectively to all the $B$-stable prime
divisors in $X$.
\end{example}

\subsection{The connected automorphism group}
\label{subsec:automorphisms}

Denote by $\Aut(X)$ the automorphism group of $X$, and by $\Aut^0(X)$
its connected component containing the identity. Since $X$ is complete
and $\Pic(X)$ is discrete, the group $\Aut^0(X)$ is linear algebraic
and its Lie algebra is the space of global vector fields, 
$\Gamma(X,\cT_X)$ (see, e.g., \cite{Ra64}). 

The subsheaf $\cS_X \subseteq \cT_X$ consisting of those vector fields
that preserve the boundary is a locally free sheaf, called the
logarithmic tangent sheaf; it contains the image of the natural map 
$\cO_X \otimes \fg \to \cT_X$, where $\fg$ denotes the Lie algebra of
$G$. The subspace $\Gamma(X,\cS_X) \subseteq \Gamma(X,\cT_X)$ is the
Lie algebra of the subgroup 
$\Aut^0(X,\partial X) \subseteq \Aut^0(X)$, and this subgroup
preserves all the $G$-orbits in $X$.

In fact, the induced map $\cO_X \otimes \fg \to \cS_X$ is surjective,
and we have an exact sequence 
\begin{equation}\label{eqn:aut}
0 \to \Gamma(X,\cS_X) \to \Gamma(X,\cT_X) \to 
\bigoplus_{i=1}^r \Gamma(X_i,\cN_{X_i/X}) \to 0
\end{equation}
where $\cN_{X_i/X} = \cO_{X_i}(X_i)$ denotes the normal sheaf. 
Moreover, the sequence
\begin{equation}\label{eqn:normal}
\CD
0 @>>> k @>{s_i}>> \Gamma(X,\cO_X(X_i)) @>>> 
\Gamma(X_i,\cN_{X_i/X}) @>>> 0
\endCD
\end{equation}
is exact for $i = 1,\ldots,r$ (for these results, see \cite{BB96}).

Together with Theorem \ref{thm:fixed}, it follows that the
$G$-module $\Gamma(X_i,\cN_{X_i/X})$ is either $0$ (if $X_i$ is
fixed), or simple with highest weight $\gamma_i$ (if $X_i$ is
base-point-free).

\begin{theorem}\label{thm:moved}
(1) The fixed divisors in $X$ are exactly the boundary divisors
$X_i$ such that $\Gamma(X,\cO_X(X_i)) = k s_i$. They are all stable
under $\Aut^0(X)$. 

\smallskip

\noindent
(2) The base-point-free boundary divisors are all moved by
$\Aut^0(X)$. For any such divisor $X_i$ and any closed
subgroup $G' \subseteq \Aut^0(X)$ containing the image of $G$, either
$X_i$ is $G'$-stable or it contains no $G'$-orbit; in particular,
$X_i$ contains no orbit of $\Aut^0(X)$. The morphism 
$\varphi : X \to X'$ associated with $X_i$ as in Theorem
\ref{thm:fixed} is equivariant under $\Aut^0(X)$, which acts
transitively on $X'$.

\smallskip

\noindent
(3) The extremal rays of $\Eff(X)_{\bbR}$ consist of the rays
generated by the fixed divisors, and of those generated 
by the pull-backs of Schubert divisors under $\Aut^0(X)$-equivariant
morphisms $\varphi : X \to X'$, where $X'$ is the quotient of
$\Aut^0(X)$ by a maximal parabolic subgroup. Both types of rays are
disjoint.
\end{theorem}

\begin{proof}
(1) Let $D$ be a fixed divisor; then $D$ is stable under $\Aut^0(X)$,
since $\xi \cdot s_D$ is a scalar multiple of $s_D$ for any 
$\xi \in \Gamma(X,\cT_X)$. In particular, $D$ is a boundary divisor;
together with Theorem \ref{thm:fixed}, this implies the assertions.

(2) Let $X_i$ be a base-point-free boundary divisor. Then $X_i$ is not
stable under $\Aut^0(X)$, by (\ref{eqn:aut}) and (\ref{eqn:normal}). 
Moreover, some positive power of $\cO_X(X_i)$ admits a linearization
for the group $\Aut^0(X)$, see \cite[Sec.~2.4]{KKLV}. Therefore, this
group acts on $X'$ so that $\varphi$ is equivariant. Let 
$G' \subseteq \Aut^0(X)$ be a closed subgroup containing the image
of $G$ and moving $X_i$. Then $G'$ acts on $X'$ and moves $X'_i$, so
that $X'$ is a unique $G'$-orbit. Thus, $X_i$ contains no $G'$-orbit. 

(3) follows from (2) together with Lemmas \ref{lem:boundary},
\ref{lem:color} and \ref{lem:mixed}.
\end{proof}

\begin{theorem}\label{thm:aut}
(1) Let $I$ be any subset of $\{1,\ldots,r\}$ containing the indices
of the fixed divisors and let $\Aut^0(X,\partial_I X)$ be the
stabilizer in $\Aut^0(X)$ of the $X_i$, $i \in I$. Then 
$\Aut^0(X,\partial_I X)$ is semi-simple and $X$ is a wonderful variety
under that group, with boundary divisors being the $X_i$, $i \in I$.

\smallskip

\noindent
(2) Let $G'$ be a closed connected subgroup of $\Aut^0(X)$ containing
the image of $G$. Then $G'$ is semi-simple and the $G'$-variety $X$ is
wonderful; its set of colors (relative to some Borel subgroup
containing the image of $B$) is exactly $\cD$. Moreover, there are
only finitely many such intermediate subgroups $G'$.

\smallskip

In particular, $\Aut^0(X)$ is semi-simple and $X$ is a wonderful
variety under that group, with boundary divisors being the fixed
divisors, and with set of colors being $\cD$.
\end{theorem}

\begin{proof} 
The Lie algebra of $\Aut^0(X,\partial_I X)$ is the kernel of the
natural map 
$$
\Gamma(X,\cT_X) \to \bigoplus_{i \in I} \Gamma(X_i,\cN_{X_i/X}).
$$
By (\ref{eqn:aut}) and Theorem \ref{thm:moved}, it follows that
$\Aut^0(X,\partial_I X)$ moves each $X_j$, $j \notin I$.  

The $G$-variety 
$$
X_I := \bigcap_{i \in I} X_i
$$ 
is wonderful and stable under $\Aut^0(X,\partial_I X)$. Moreover, any
proper $G$-stable closed subvariety of $X_I$ is contained in $X_j$ for
some $j \notin I$, and hence contains no orbit of 
$\Aut^0(X,\partial_I X)$, by Theorem \ref{thm:moved} again. It
follows that $X_I$ is a unique orbit of $\Aut^0(X,\partial_I X)$. 

Let $G_I$ be a Levi subgroup of $\Aut^0(X,\partial_I X)$ containing
the image of $G$ and let $B_I$ be a Borel subgroup of $G_I$ containing
the image of $B$. Then $G_I$ is semi-simple (since the
$G$-automorphism group of $X$ is finite) and acts transitively on the
projective $\Aut^0(X,\partial_I X)$-orbit $X_I$. Clearly, the
$G_I$-variety $X$ is spherical, $X_I$ is its unique closed orbit, and
its $B_I$-stable prime divisors are among the $D \in \cD$ and
$X_1,\ldots,X_r$. By Theorem \ref{thm:moved}, the boundary divisors
are either $G_I$-stable or contain no $G_I$-orbit. Thus, the
$G_I$-variety $X$ is toroidal, and hence wonderful by
\cite[Sec.~7]{Kn91}. Its boundary divisors are those $X_i$ that contain
$X_I$, i.e., such that $i \in I$. They are all stable under
$\Aut^0(X,\partial_I X)$, so that this group is semi-simple by Lemma 
\ref{lem:ss} below. In other words, 
$G_I = \Aut^0(X,\partial_I X)$. This completes the proof of (1).

For (2), let $I$ be the set of $G'$-stable boundary divisors. Arguing
as above, we see that $X_I$ is a unique orbit of $G'$, and hence of a
Levi subgroup $G''\subseteq G'$ containing the image of $G$. Moreover,
$G''$ is semi-simple and $X$ is wonderful for that group, with
boundary divisors indexed by $I$. Applying Lemma \ref{lem:ss} again,
we obtain that $G'' = G'$. The finiteness of the set of subgroups $G'$
follows from Lemma \ref{lem:finite} below.

To complete the proof, it suffices to show that $X_B = X_{B'}$ for an
appropriate Borel subgroup $B'$ of $G'$ containing the image of
$B$. By \cite[Sec.~1]{Lu97}, $X_B$ is the Bialynicki--Birula cell 
$X(\lambda,z)$ where $z \in X_B$ is the unique $T$-fixed point, and
$\lambda$ is a general one-parameter subgroup of $T$ such that the
associated parabolic subgroup $G(\lambda)$ equals $B$.  This defines a
parabolic subgroup $G'(\lambda)$ of $G'$, with Levi subgroup being the
centralizer $C_{G'}(\lambda) = C_{G'}(T)$. 

We claim that $X_B$ is stable under $G'(\lambda)$. Given $x \in X_B$
and $g' \in G'(\lambda)$,  
$$
\lambda(t) \, g' \cdot x \to g'_0 \cdot z \quad \text{as} \quad t \to 0,
$$
where $g'_0 := \lim_{t \to 0} \lambda(t) g' \lambda(t^{-1})$ is in
$C_{G'}(\lambda)$. Moreover, $C_{G'}(\lambda)$ is connected, and 
$C_{G'}(\lambda) \cdot z$ consists of $T$-fixed points. The latter
being isolated, it follows that $g'_0 \cdot z = z$ which proves the
claim. 

Let $B'$ be a Borel subgroup of $G'(\lambda)$ containing the image of
$B$. Then $X_B$ is stable under $B'$, which implies the desired equality.
\end{proof}

\begin{lemma}\label{lem:ss}
Any closed connected subgroup of $\Aut^0(X, \partial X)$ which
contains $G$ is semi-simple.
\end{lemma}

\begin{proof}
Otherwise, there exists a non-trivial closed connected unipotent
subgroup $A$ of $\Aut^0(X,\partial X)$ which is normalized by
$G$; we may further assume that $A$ is abelian. The semi-direct
product $\tG := G \cdot A$ acts on $X$ with the same orbits as those
of $G$. Moreover, the invariant subgroup $A^G$ is trivial, since
$\Aut_G(X)$ is finite.

Let $\tfg = \fg \oplus \fa$ be the Lie algebra of $\tG$ and
let $S_X$ be the $\tG$-linearized vector bundle on $X$ with sheaf of
local sections $\cS_X$. Then the map $\fg \to \Gamma(X,\cS_X)$ yields
a $G$-equivariant morphism
$$
\mu: S_X^* \to \fg^*,
$$
the localized moment map of \cite{Kn94}. This morphism is proper, and 
so is the analogous $\tG$-equivariant morphism
$$
\tmu: S_X^* \to \tfg^*.
$$
Thus, both images $\Ima(\mu)$, $\Ima(\tmu)$ are closed, and the
projection 
$$
p : \tfg^* = \fg^* \oplus \fa^* \to \fg^*
$$
restricts to a finite $G$-morphism, also denoted by
$p : \Ima(\tmu) \to \Ima(\mu)$.

Consider the quotient map 
$$
q : \fg^* \to \fg^*/\!/G.
$$ 
We claim that the composite morphism 
$$
q \circ \mu : S_X^* \to \fg^* \to \fg^*/\!/G
$$
is $\tG$-invariant. Indeed, by \cite[Thm.~6.1]{Kn94}, any regular
function on $S_X^*$ arising from $q \circ \mu$ is the symbol of a
$G$-invariant differential operator on $X$ which is completely regular
in the sense of [loc.~cit., p.262]. Thus, this operator commutes with the
action of $\tfg \subseteq \Gamma(X,\cS_X)$, by 
[loc.~cit., Cor.7.6]. This implies our claim. 

By that claim, the composite morphism 
$$
\pi:= q \circ p : \Ima(\tmu) \to \Ima(\mu) \to \fg^*/\!/G
$$
is $\tG$-invariant. Moreover, since $p$ is finite and any fiber of $q$ 
contains only finitely many $G$-orbits, any fiber of $\pi$
contains only finitely many $G$-orbits as well. It follows that each
general $\tG$-orbit in $\Ima(\tmu)$ contains an open $G$-orbit. In other
words, the equality $\tfg \cdot f = \fg \cdot f$ holds for any $f$ in
a dense subset of $\Ima(\tmu)$. By \cite[Thm.~2.2]{BK94}, this implies
that $f \in \fg^*$, i.e., $\Ima(\tmu) \subseteq \fg^*$. In other
words, the map $\fa \to \Gamma(X,\cS_X)$ is zero, a contradiction.
\end{proof}

\begin{lemma}\label{lem:finite}
Let $G_1 \subseteq G_2$ be semi-simple groups such that any closed
connected subgroup $G' \subseteq G_2$ containing $G_1$ is
semi-simple. Then the number of such intermediate subgroups is 
finite. 
\end{lemma}

\begin{proof}
This result is certainly known; for lack of a reference, we provide an
argument. Consider the semi-simple Lie algebras 
$\fg_1 \subseteq \fg' \subseteq \fg_2$. By \cite{Ri67}, the orbit 
$G_2 \cdot \fg'$ is open in the variety of Lie subalgebras of $\fg_2$,
and hence the number of $G_2$-conjugacy classes of intermediate
subgroups $G'$ is finite. Moreover, the centralizer of $G_1$ in
$G_2$ is finite. By \cite{Ri82}, it follows that any $G_2$-orbit
contains only finitely many $G_1$-fixed points. Hence any conjugacy
class contains only finitely many overgroups of $G_1$.
\end{proof}

\begin{example}\label{ex:group4}
We determine the automorphism group of $X = \bar{G}_{\ad}$. Write 
$$
G = (\SL_2)^n \times G',
$$ 
where $G'$ contains no direct factor isomorphic to $\SL_2$. Then 
$$
X \cong \bbP(M_2)^n \times X' \cong (\bbP^3)^n \times X',
$$
where $X' := \bar{G'}_{\ad}$. Thus,
$$
\Aut^0(X) \cong (\PSL_4)^n \times \Aut^0(X').
$$
Since all boundary divisors of $X'$ are fixed, the semi-simple group
$\Aut^0(X')$ stabilizes all the $G' \times G'$-orbits; as a
consequence, it acts faithfully on the closed orbit. The latter is
isomorphic to $G'/B' \times G'/B'$, and
$\Aut^0(G'/B' \times G'/B') = G'_{\ad}\times G'_{\ad}$ by \cite{De77},
so that
$$
\Aut^0(X') = G'_{\ad} \times G'_{\ad}
$$
as well.

The full automorphism group $\Aut(X)$ acts on the nef monoid and
permutes its generators $[D]$, $D\in \cD$; it stabilizes the subset of
those $[D]$ that generate an extremal ray of $\Eff(X)_{\bbR}$, i.e.,
the pull-backs of hyperplanes under the $n$ projections
$X \to \bbP^3$. It follows that 
$$
\Aut(X) \cong \Aut((\bbP^3)^n) \times \Aut(X').
$$
Moreover, $\Aut((\bbP^3)^n)$ is the semi-direct product of $(\PSL_4)^n$
by the symmetric group $S_n$ permuting the copies of $\bbP^3$. 

To determine $\Aut(X')$, we consider its action by conjugation on its
connected component $ G'_{\ad} \times G'_{\ad}$. This action is
faithful, since $\Aut_{G' \times G'}(X')$ is trivial. Moreover, any  
$\varphi \in \Aut(X')$ sends $\diag(G'_{\ad})$ to a conjugate subgroup
in $G'_{\ad} \times G'_{\ad}$, since $\varphi$ preserves the open
orbit 
$X'^0 \cong  (G'_{\ad} \times G'_{\ad})/\diag(G'_{\ad})$. It follows
that the group $\Aut(X')$ is generated by $G'_{\ad} \times G'_{\ad}$
together with the flip $(g,h) \mapsto (h,g)$ and with the
automorphisms $(\varphi,\varphi)$, where $\varphi$ arises from an
automorphism of the Dynkin diagram of $G'$.
\end{example}

\section{The total coordinate ring of a wonderful variety}
\label{sec:wonderful}
 
\subsection{Definitions and basic properties}
\label{subsec:basic}

The notation and assumptions are as in Section \ref{sec:divisors}; in
particular, $X$ denotes a wonderful variety under the simply-connected
semi-simple group $G$, and $\cD$ denotes its set of colors.

Consider the sum
$$
\bigoplus_{(n_D) \in \bbZ^{\cD}} 
\cO_X (\sum_{D \in \cD} n_D D ).
$$
This is a sheaf of $\cO_X$-algebras, graded by $\bbZ^{\cD}$, and
compatibly $G$-linearized. It defines a scheme $\hX$ equipped
with a morphism
$$
p: \hX \to X
$$
and with an action of the group 
$$
\tG := G \times \bbG_m^{\cD}
$$
where $\bbG_m^{\cD} \cong \Hom(\Pic(X), \bbG_m)$ denotes the torus with
character group $\bbZ^{\cD} \cong \Pic(X)$. The morphism $p$ is
$G$-equivariant and is a principal $\bbG_m^{\cD}$-bundle: the 
\emph{universal torsor} over $X$. In particular, $\hX$ is a
non-singular spherical variety under the connected reductive group
$\tG$.

By \cite[Prop.~3.10]{BH03}, the variety $\hX$ is quasi-affine. 
This may also be seen directly: let $\cL$ be a very ample invertible
sheaf on $X$. Then the sheaves $\cL \otimes \cO_X(D)$, $D \in \cD$ are
all very ample, and form a basis of $\Pic(X)$. Consider the
corresponding projective embeddings $X \hookrightarrow \bbP_D$, and
the tautological principal $\bbG_m$-bundles 
$\widehat{\bbP}_D \to \bbP_D$. Then each $\widehat{\bbP}_D$ is the
complement of a point in an affine space, and $\hX$ is the
pull-back of the product bundle under the diagonal embedding 
$X \to \prod_{D \in \cD} \bbP_D$.

Next let
$$
R(X) := \Gamma(\hX, \cO_{\hX}) =
\bigoplus_{(n_D) \in \bbZ^{\cD}} 
\Gamma ( X,\cO_X(\sum_{D \in \cD} n_D D)).
$$
This is a $\bbZ^{\cD}$-graded $k$-algebra called the 
\emph{total coordinate ring} of $X$. The set of degrees of its
homogeneous elements is the monoid $\Eff(X)$, which generates the
rational polyhedral cone $\Eff(X)_{\bbR}$ containing no line. As a
consequence, $R(X)$ admits a coarser grading by the non-negative
integers, such that $R(X)_0 = \Gamma(X,\cO_X) = k$. It follows that
$R(X)^* = k^*$.

The action of $\tG$ on $\hX$ yields an action on $R(X)$. 
Since $\hX$ is spherical, the $k$-algebra $R(X)$ is finitely
generated (see \cite{Kn93}) and normal. Thus,
$$
\tX := \Spec R(X)
$$
is a normal affine variety equipped with a $\tG$-equivariant morphism 
$$
\iota : \hX \to \tX.
$$
By \cite{BH03}, $\iota$ is an open immersion. Again, this may be
seen directly: let $\cL$ be as above, then $\hX$ is covered by
the affine open subsets $\hX_s := p^{-1}(X_s)$, where 
$s$ is a non-zero global section of $\cL$. Moreover, 
$\hX_s = \iota^{-1}(\tX_s)$ where $s$ is regarded as a
homogeneous element of $R(X)$, and the induced homomorphism
$R(X)[s^{-1}] \to k[\hX_s] = k[p^{-1}(X_s)]$ 
is an isomorphism, since each space 
$\Gamma(X_s, \cO_X(\sum n_D D))$ is the union of its subspaces
$\Gamma(X, \cL^N(\sum n_D D)) s^{-N}$, where $N \ge 0$.

Since $\hX$ and $\tX$ have the same algebra of regular functions,
the closed subset $\tX \setminus \iota(\hX)$ has codimension at least
$2$ in $\tX$. Thus, $\tX$ is a spherical $\tG$-variety with open
$\tG$-orbit
$$
\tX^0 := p^{-1}(X^0).
$$ 
Moreover, 
$$
\tX^0_{\tB} := p^{-1}(X^0_B)
$$ 
is the open orbit of 
$$ 
\tB := B \times \bbG_m^{\cD}
$$
(a Borel subgroup of $\tG$), and 
\begin{equation}\label{eqn:tdiv}
\tX \setminus \tX^0_{\tB}  = 
\tX_1 \cup \cdots \cup \tX_r \cup \bigcup_{D \in \cD} \tD,
\end{equation}
where $\tX_i$ denotes the closure in $\tX$ of the prime divisor 
$p^{-1}(X_i)$, and likewise for $\tD$. Note also that $\tX$ contains a
fixed point of $\tG$, since $R(X)$ admits a positive $\tG$-invariant
grading. 

By \cite[Prop.~8.4]{BH03} or \cite[Cor.~1.2]{EKW04}, $R(X)$ is a
UFD. This may also be seen directly, by showing that the divisor
class group of $\tX$ is trivial. Indeed, the open subset $\tX^0_{\tB}$
has a trivial divisor class group, and the ideal of each irreducible
component of its complement is generated by the corresponding
canonical section. 

We now describe the invariant subalgebras $R(X)^U$ and $R(X)^G$. 

\begin{proposition}\label{prop:quot}
(1) The algebra $R(X)^U$ is freely generated by the canonical sections
$s_1,\ldots,s_r$ and the $s_D$, $D \in \cD$. Each $s_D$ is an
eigenvector of $\bbG_m^{\cD}$ of weight the $D$-component 
$$
\varepsilon_D : \bbG_m^{\cD} \to \bbG_m.
$$
Moreover, each $s_i$ is a $\bbG_m^{\cD}$-eigenvector of weight
$\sum_{D \in \cD} \langle D,\gamma_i \rangle \, \varepsilon_D$.

\smallskip

\noindent
(2) The subalgebra $R(X)^G$ is freely generated by $s_1,\ldots,s_r$.

\smallskip

\noindent
(3) $R(X)$ is a free module over $R(X)^G$.

\smallskip

\noindent
(4) The quotient morphism
$$
q : \tilde X = \Spec R(X) \to \Spec R(X)^G =: \tilde X /\!/G \cong \bbA^r
$$
is flat and its (scheme-theoretic) fibers are normal varieties.
\end{proposition}

\begin{proof}
(1) Any element of $R(X)^U$ can be uniquely decomposed into a sum of
$\tilde B$-eigenvectors. We claim that any such eigenvector $\sigma$
is a monomial in $s_1,\ldots,s_r$ and the $s_D$, $D \in \cD$, with
uniquely determined exponents. Indeed, $\sigma$ is a $B$-eigenvector
in some space $\Gamma(X,\cO_X(\sum_{D \in \cD} n_D D))$. Its divisor
of zeroes may be written uniquely as
$$
\div (\sigma) = \sum_{D \in \cD} \ord_D(\sigma) \, D + 
\sum_{i=1}^r \ord_{X_i}(\sigma) \, X_i.
$$
Thus, $\sigma$ is a scalar multiple of the monomial 
$\prod_{D \in \cD} s_D^{\ord_D(\sigma)} \cdot
\prod_{i=1}^r s_i^{\ord_{X_i}(\sigma)}$
and of no other monomial. This proves the claim and, in turn, the
assertion on $R(X)^U$. 

Clearly, each $s_D$ has weight $\varepsilon_D$; by (\ref{eqn:irel}),
this implies the assertion on the weight of $s_i$.

(2) Likewise, each $\tG$-eigenvector in $R(X)$ is a monomial in
$s_1,\ldots,s_r$ with uniquely determined exponents.

(3) Let $M$ be the $G$-submodule of $R(X)$ generated by the
monomials in the $s_D$, $D \in \cD$. Then $M^U$ is just the
polynomial ring $k[s_D, D\in \cD]$ (but $M$ is generally not a subring 
of $R(X)$). Consider the multiplication map
$$
m : R(X)^G \otimes M \to R(X).
$$
This is a morphism of $G$-modules, which restricts to an isomorphism
\begin{equation}\label{eqn:Uinv}
(R(X)^G \otimes M)^U = R(X)^G \otimes M^U \cong R(X)^U
\end{equation}
by (1). Thus, $m$ is an isomorphism.

(4) By (3), $q$ is flat. Let $F$ be a fiber at a closed point of
$\tilde X/\!/G$, then $F$ is an affine $G$-scheme and the restriction 
map $M^U \to k[F]^U$ is an isomorphism by (\ref{eqn:Uinv}). Thus,
$k[F]^U$ is a polynomial ring. It follows that $F$ is normal; see,
e.g., \cite[Sec.~6]{Po87}.
\end{proof}

The base point $x \in X^0$ yields a base point 
$\tx \in \tX^0 = \hX^0$, as follows : $\hX^0$ is the 
total space of the direct sum of the line bundles on $X^0 = G/H$
associated with the invertible sheaves $\cO_{G/H}(D)$, $D \in \cD$,
minus the union of the partial sums. Moreover, we may regard each
$f_D$ as a global section of $\cO_{G/H}(D)$ which does not vanish at
$x$ (since $f_D(1) = 1$). Thus, there is a unique 
$\tx \in \hX^0$ such that $f_D(\tx) = 1$ for all $D$. 

Clearly, the isotropy group $\tH := \tG_{\tx}$ is the image of the
homomorphism
$$
H \to G \times \bbG_m^{\cD}, \quad 
h \mapsto (h, (\chi_D(h), D \in \cD)). 
$$
The combinatorial data of the spherical homogeneous space $\tG/\tH$
(that is, its weight group, colors, and valuation cone) are described
in the following statement, which follows easily from Proposition 
\ref{prop:quot}.

\begin{lemma}\label{lem:comb}
(1) The pairs $(\omega_D,\varepsilon_D)$, $D \in \cD$, and the pairs
$(0,\sum_D \langle D,\gamma_i\rangle \, \varepsilon_D)$,
$i=1,\ldots,r$, form a basis of the weight group 
$\tLambda(\tX) \subseteq \tLambda := \Lambda \times \bbZ^{\cD}$.
The dual basis of $\Hom(\tLambda(\tX),\bbZ)$ consists of the 
$\trho(v_{\tD})$, $D \in \cD$, and $\tv_1,\ldots,\tv_r$. The map 
$\trho: \tcD = \cD \to \Hom(\tLambda(\tX),\bbZ)$ is injective.

\smallskip

\noindent
(2) The projection $p : \hX \to X$ yields a morphism 
$p^*: \Lambda(X) \to \tLambda(\hX) = \tLambda(\tX)$ which sends
any spherical root $\gamma_i$ to 
$$
(\gamma_i,0) = 
\sum_{D \in \cD} \langle D,\gamma_i \rangle \, (\omega_D,\varepsilon_D) 
- (0, \sum_{D \in \cD} \langle D,\gamma_i \rangle \, \varepsilon_D).
$$ 
The $(\gamma_i,0)$ and the $(\omega_D,\varepsilon_D)$, $D \in \cD$,
form another basis of $\Hom(\tLambda(\tX),\bbZ)$.

\smallskip

\noindent
(3) The transpose
$p_*: \Hom(\tLambda(\tX),\bbR) \to \Hom(\Lambda(X),\bbR)$
maps each $\trho(v_{\tD})$ to $\rho(v_D)$, and each $\tv_i$ to
$v_i$. Moreover, the valuation cone $\tcV$ is the preimage of
$\cV$ under $p_*$.
\end{lemma}

Together with the classification of embeddings of spherical
homogeneous spaces (see \cite[Sec.~4]{Kn91}), this description implies
readily the following:

\begin{lemma}\label{lem:orbits}
The $\tG$-orbits in $\tX$ correspond bijectively to those pairs 
$(\cE,I)$, where $\cE \subseteq \cD$, $I \subseteq \{ 1,\ldots, r\}$,
such that the relative interior of  
$\Cone(\rho(v_D), D\in\cE ; v_i, i \in I)$
meets the valuation cone $\cV$. Equivalently, there exist positive
integers $x_D$, $D \in \cE$, such that 
$\sum_{D \in \cE} \, x_D \, \langle D,\gamma_i \rangle \leq 0$ for all 
$i \notin I$.

Under this correspondence $(\cE,I) \leftrightarrow \cO_{\cE,I}$, the
ideal of the closure $\bar{\cO}_{\cE,I}$ is generated as a $G$-module
by the monomials 
$\prod_{i=1}^r s_i^{n_i} \prod_{D \in \cD} s_D^{n_D}$, where $n_i >0$
for all $i \in I$, and $n_D >0$ for all $D \in \cE$. 

In particular, $q(\bar{\cO}_{\cE,I})$ is the linear subspace of
$\bbA^r$ defined by the vanishing of the $s_i$, $i \in I$. Moreover,
the orbit $\cO_{\cF,J}$ is contained in $\bar{\cO}_{\cE,I}$ if and
only if $\cE \subseteq \cF$ and $I \subseteq J$. The orbits contained
in $\hX$ correspond to the pairs $(\emptyset,I)$, where $I$ is an
arbitrary subset of $\{ 1,\ldots, r\}$.
\end{lemma}

\begin{example}\label{ex:group5}
If $X = \bar{G}_{\ad}$, then $\tX$ is nothing but the enveloping
semigroup $\Env(G)$ of Vinberg \cite{Vi95b}, as proved in
\cite[Thm.~20]{Ri97} (see also \cite{Ri02}; another proof will be given
in Example \ref{ex:group6}). Moreover, the quotient morphism 
$q: \tX \to \bbA^r$ is the abelianization of $\Env(G)$ defined in
\cite[p.~149]{Vi95b}. Lemma \ref{lem:orbits} gives back Vinberg's
description \cite[Thm.~6]{Vi95b} of the $G\times G$-orbits in
$\Env(G)$, see \cite[Sec.~5.3]{Ri97} for details.
\end{example}

\subsection{An algebraic description}
\label{subsec:algebraic}

Consider the action of $G \times T$ on $\tX$, where $G$ acts
naturally and $T$ acts via the homomorphism
$$
u : T \to \bbG_m^{\cD}, \quad t \mapsto (\omega_D(t))_{D\in \cD}.
$$
The open $G \times \bbG_m^{\cD}$-orbit $\tX^0$ is stable under
this action; its structure is described by the following

\begin{lemma}\label{lem:fib}
There is an isomorphism of $G \times T$-varieties 
$$
\tX^0 \cong G/K \times^{T\cap H} T
$$
where $K$ denotes the intersection of the kernels of all characters of
$H$ (this is a normal subgroup of $H$), and $T\cap H$ acts on $G/K$
via the restriction of the right action of $H$ on $G/K$.

Via this isomorphism, the restriction of the quotient map 
$q^0: \tX^0 \to \bbA^r$ is identified with the projection 
$G/K \times^{T\cap H} T \to T/T\cap H$, where $T/T\cap H$ is
isomorphic to $\bbG_m^r \subset \bbA^r$ via the spherical roots
$\gamma_1,\ldots,\gamma_r$.
\end{lemma}

\begin{proof}
We have an isomorphism of $G\times \bbG_m^{\cD}$-varieties 
$$
\tX^0 \cong G \times ^H \bbG_m^{\cD}
$$ 
where $H$ acts on $\bbG_m^{\cD}$ via its characters $\chi_D$. This
identifies $q^0$ to the $H$-invariant map 
$$
G \times \bbG_m^{\cD} \to \bbA^r, \quad
(g;t_D, D\in \cD) \mapsto 
(\prod_{D \in \cD} t_D^{\langle D,\gamma_i \rangle})_{i=1,\ldots,r}
$$
with image $\bbG_m^r$. Since the subgroup $K\subseteq H$ acts trivially
on $\bbG_m^{\cD}$, we obtain 
$$
\tX^0 \cong G/K \times^{H/K} \bbG_m^{\cD}.
$$ 
Moreover, we have a commutative diagram
$$
\CD
1 @>>> T \cap H @>>> T           @>>> T/T \cap H @>>> 1 \\
& &       @VVV       @V{u}VV          @V{\gamma_1,\ldots,\gamma_r}VV \\  
1 @>>> H/K      @>>> \bbG_m^{\cD} @>>> \bbG_m^r    @>>> 1  \\
\endCD
$$
The square on the left is a fiber square, since all involved groups
are diagonalizable and the dual square of character groups is
$$
\CD
\bbZ^{\cD} @>>> \cX(H/K) = \cX(H)                     \\
@VVV           @VVV                                  \\
\cX(T) = \cX(B)    @>>> \cX(T \cap H) = \cX(B\cap H) \\ 
\endCD
$$
which is cartesian by Lemma \ref{lem:str}.
It follows that the natural map
$$
G/K \times^{T \cap H} T \to G/K \times^{H/K} \bbG_m^{\cD}
$$
is an isomorphism.
\end{proof}

The subgroup $K$ plays an important role in \cite[Sec.~6]{Lu01}. As
already observed, the quotient $H/K$ is diagonalizable with character
group $\cX(H)$. Moreover, the homogeneous space $G/K$ is quasi-affine
and satisfies
$$
k[G/K] = \bigoplus_{\chi \in \cX(H)} k[G]^{(H)}_{\chi}.
$$
In particular, the $G \times (H/K)$-module $k[G/K]$ is
multiplicity-free, and the algebra $k[G/K]^U = k[G]^{U\times K}$ is a 
polynomial ring in the $f_D$, $D \in \cD$. It follows that the
algebra $k[G/K]$ is finitely generated; see, e.g., \cite[Cor.~4]{Po87}.

We now derive from Lemma \ref{lem:fib} an algebraic description of the
coordinate ring of the quasi-affine variety $\tX^0$. For this, we
introduce some notation: for any $G$-module $M$ and any dominant
weight $\lambda$, we denote by $M_{(\lambda)}$ the isotypical
component of $M$ of type $\lambda$, that is, the sum of all simple
$G$-submodules with highest weight $\lambda$. Also, we denote by
$(e^{\mu}, \mu \in \Lambda)$ the basis of the $k$-vector space $k[T]$
consisting of characters.

\begin{lemma}\label{lem:alg}
There is a isomorphism of $G \times T$-algebras
\begin{equation}\label{eqn:is0}
k[\tX^0] \cong 
\bigoplus_{\lambda \in \Lambda^+, \; \mu \in \Lambda, \; 
\lambda - \mu \in \Lambda(X)} 
k[G/K]_{(\lambda)} \, e^{\mu},
\end{equation}
where the right-hand side is a subalgebra of 
$k[G \times T] = \bigoplus_{\mu \in \Lambda} k[G] \, e^{\mu}$.
\end{lemma}

\begin{proof}
From Lemma \ref{lem:fib} we obtain
$$
k[\tX^0] \cong (k[G/K] \otimes k[T])^{T\cap H}
\cong \bigoplus k[G]^{(H)}_{\chi} \, e^{\mu},
$$
where the sum runs over those pairs 
$(\chi,\mu) \in \cX(H) \times \Lambda$ such that 
$\chi \vert_{T\cap H} = \mu \vert_{T\cap H}$.

Moreover, given $\nu \in \cX(T\cap H)$, the equality
$$
\bigoplus_{\chi} k[G]^{(H)}_{\chi} =
\bigoplus_{\lambda} k[G]^{(H)}_{(\lambda)}
$$
holds, where the sum on the left-hand side (resp.~right-hand side)
runs over those $\chi \in \cX(H)$ (resp.~$\lambda \in \Lambda^+$)
such that $\chi\vert_{T\cap H} = \nu$ 
(resp. $\lambda \vert_{T\cap H} = \nu$). Indeed,
both sides are $G$-submodules of $k[G/K]$; and by Lemma \ref{lem:str},
they have the same $B$-eigenvectors, namely, those monomials in the
$f_D$, $D \in \cD$, that have weight $\nu$ with respect to $T\cap H$.

Finally, the condition that 
$\lambda\vert_{T\cap H} = \mu\vert_{T\cap H}$ 
is equivalent to $\lambda - \mu \in \Lambda(X)$, since 
$\Lambda(X) = \cX(T)^{T\cap H}$. Putting these facts together yields 
the desired isomorphism. 
\end{proof}

Next we deduce from Lemma \ref{lem:alg} an algebraic description of
the ring $R(X)$; to formulate it, we need an additional notation. 
Given $\lambda, \mu \in \Lambda$, we write
$$
\lambda \leq_X \mu
$$ 
if the difference $\mu - \lambda$ is a linear combination of spherical 
roots with non-negative coefficients; then 
$\mu - \lambda \in \Lambda(X)$, and $\lambda \leq \mu$ for the usual
ordering on weights.
 
\begin{theorem}\label{thm:rees}
There is an isomorphism of $G \times T$-algebras
\begin{equation}\label{eqn:is}
R(X) \cong 
\bigoplus_{\lambda \in \Lambda^+, \; \mu \in \Lambda, \; \lambda \leq_X \mu}
 k[G/K]_{(\lambda)} \, e^{\mu},
\end{equation}
where the left-hand side is a subalgebra of $k[G \times T]$.
This isomorphism identifies each $s_i$ with $e^{\gamma_i}$ and each
$s_D$ with $f_D \, e^{\omega_D}$.
\end{theorem}

\begin{proof}
Via the isomorphism (\ref{eqn:is0}), the restriction
$s_i\vert_{\tX^0}$ is identified with $e^{\gamma_i}$, since both span
the subspace of $G \times T$-eigenvectors of weight
$(0,\gamma_i)$, and both take the value $1$ at the base point
$\tx$. Likewise, the restriction $s_D\vert_{\tX^0}$ is identified with
$f_D \, e^{\omega_D}$. This implies (\ref{eqn:is}), since $R(X)$
is the $G$-submodule of $k[\tX^0]$ generated by the monomials in 
$s_1, \ldots, s_r$ and the $s_D$, $D\in \cD$.  
\end{proof}

Since $R(X)$ is a subalgebra of $k[G \times T]$, we must have
$$
k[G/K]_{(\lambda)} \; k[G/K]_{(\mu)} \subseteq 
\bigoplus_{\nu \in \Lambda^+,\, \nu \leq_X \lambda + \mu}
k[G/K]_{(\nu)}
$$
for any dominant weights $\lambda,\mu$, where the left-hand side
denotes the product in $k[G/K]$. In other words, the 
\emph{root monoid} of the algebra $k[G/K]$ (that is, the submonoid of
$\Lambda$ generated by the differences $\lambda + \mu - \nu$, where 
$\lambda, \mu, \nu$ are dominant weights such that 
$k[G/K]_{(\lambda)} \, k[G/K]_{(\mu)}$ has a non-zero projection to
$k[G/K]_{(\nu)}$) is contained in the monoid generated by the
spherical roots. In fact, both monoids generate the same cone in
$\Lambda(X)_{\bbR}$ (see, e.g., \cite[Lem.~6.1]{Kn91}; this result will
be strengthened in Proposition \ref{prop:gs}). Thus, the subspaces 
$$
k[G/K]_{\leq_X \mu} := 
\bigoplus_{\lambda \in \Lambda, \, \lambda \leq_X \mu}
k[G/K]_{(\lambda)}, \quad \mu \in \Lambda
$$
form an ascending filtration of $k[G/K]$ indexed by the partially
ordered set $\Lambda$, and the right-hand side of (\ref{eqn:is}) is
the associated Rees algebra.

Theorem \ref{thm:rees} yields a description of the fibers of the
quotient map $q : \tX \to \bbA^r$ at all closed points. Indeed, let
$\tX_x$ be the fiber at $x = (x_1,\ldots,x_r) \in \bbA^r$ and let $I =
I(x)$ be the set of indices $i$ such that $x_i \neq 0$. Then, by
(\ref{eqn:is}), the algebra
$k[\tX_x] \cong R(X)/(s_1 - x_1,\ldots,s_r - x_r)$
is isomorphic to $k[G/K]$ endowed with the new multiplication
$$
m_I : k[G/K] \otimes k[G/K] \to k[G/K]
$$ 
such that the restriction of $m_I$ to any 
$k[G/K]_{(\lambda)} \otimes k[G/K]_{(\mu)}$
is the original multiplication 
$$
m: k[G/K]_{(\lambda)} \otimes k[G/K]_{(\mu)} \to
\bigoplus_{\nu \in \Lambda^+, \, \nu \leq_X \lambda + \mu}
k[G/K]_{(\nu)}
$$
followed by the projection onto the partial sum of those
$k[G/K]_{(\nu)}$ such that $\lambda + \mu - \nu$ is a linear
combination of the $\gamma_i$, $i \in I$.

In particular, the general fibers $\tX_x$, $x \in \bbG_m^r$, satisfy 
$k[\tX_x] \cong k[G/K]$ as $G$-algebras. In other words,
these fibers are all isomorphic to the \emph{canonical embedding}
$$
\CE(G/K) := \Spec k[G/K].
$$  
This is a normal affine $G$-variety in which the homogeneous space $G/K$
is embedded as an open orbit with complement of codimension $\geq 2$
(since $G/K$ is quasi-affine and its coordinate ring is finitely
generated). 

On the other hand, in the coordinate ring of the special fiber
$\tX_0$, the product of any two simple $G$-submodules
$V(\lambda),V(\mu)$ is just their Cartan product 
$V(\lambda + \mu)$. Thus, the $G$-variety $\tX_0$ is 
\emph{horospherical}, i.e., the $G$-isotropy group of any point
contains a maximal unipotent subgroup; see \cite[Sec.~4]{Po87}. In
fact, the fibers of $q$ realize a degeneration of $\CE(G/K)$ to its
``horospherical contraction'' studied in [loc.~cit.]. Also, note that
$\tX_0$ is the unique horospherical fiber.

\begin{example}\label{ex:group6}
If $X= \bar{G}_{\ad}$, then we have to replace the acting group $G$
with $G\times G$, and the subgroup $K$ with $\diag(G)$. Then the ring
$k[G/K]$ is replaced with $k[(G\times G)/\diag(G)] = k[G]$, which is
isomorphic to $\bigoplus_{\lambda \in \Lambda} \End(V(\lambda))$
as a $G\times G$-module. It follows that
$$
R(X) \cong 
\bigoplus_{\lambda \in \Lambda^+,\,\mu \in \Lambda, \, \lambda \leq \mu}
\End(V(\lambda)) \, e^{\mu},
$$
where $\leq$ is the usual ordering on weights. By \cite[p.~152]{Vi95b},
this algebra is the coordinate ring of the enveloping monoid
$\Env(G)$. The horospherical degeneration of $\tX$ is the asymptotic
semigroup of \cite{Vi95a}; the other fibers of the quotient morphism
$q$ are described in \cite[Lem.~7.17]{AB04}.

By \cite{Vi95b} again, the unit group of $\Env(G)$ is $G \times^Z T$,
where the center $Z$ acts diagonally on $G \times T$. Alternatively,
this follows from Lemma \ref{lem:fib}, which also implies that the
action of $G \times T$ on $R(X)$ (where $X$ is any wonderful
$G$-variety) factors through an action of $G \times^Z T$. In fact, the
latter action extends to an action of the monoid $\Env(G)$ (as follows
from Theorem \ref{thm:rees}) which may deserve a detailed investigation. 
\end{example}

\subsection{Generators and relations}
\label{subsec:relations}

For any $D \in \cD$, let $V_D \cong V(\omega_D)$ be the $G$-submodule  
of $\Gamma(X,\cO_X(D))$ generated by $s_D$. Proposition
\ref{prop:quot} implies readily that the algebra $R(X)$ is generated
by $s_1, \ldots, s_r$ and the $G$-modules $V_D$, where $D\in \cD$.
Thus, we may write $R(X) = S/I$, where $I$ is an ideal of 
$$
S := k[s_1, \ldots, s_r] \otimes \Sym(\bigoplus_{D \in \cD} V_D).
$$
To construct generators of $I$, we determine the product in $R(X)$ of
any two submodules $V_{D_1}, V_{D_2}$, as a quotient of their tensor
product $V_{D_1}\otimes V_{D_2}$.

Regard each $V_D$ as the submodule of $k[G/K]$ generated by $f_D$. More
generally, given a family of non-negative integers 
$(n_D)_{D \in \cD}$, let $V_{\sum n_D D}$ be the $G$-submodule of
$k[G/K]$ generated by $\prod f_D^{n_D}$. Then 
$$
k[G/K] = \bigoplus_{(n_D)} V_{\sum n_D D},
$$
where the sum runs over all families of non-negative integers.
Thus, any product $V_{D_1} \, V_{D_2}$ in $k[G/K]$ may be decomposed
into a partial sum of $V_{\sum n_D D}$'s, where
$$
\sum_{D \in \cD} n_D \, \omega_D \leq_X \omega_{D_1} + \omega_{D_2}
\quad \text{and} \quad
\sum_{D \in \cD} n_D \, \chi_D = \chi_{D_1} + \chi_{D_2}.
$$
In particular, we may index the components $V_{\sum n_D D}$ of 
$V_{D_1} \, V_{D_2}$ by their highest weight 
$\lambda = \sum_{D \in \cD} n_D \, \omega_D$. 
Moreover, each component $V(\lambda)$ may be embedded into the
tensor product $\bigotimes_{D\in \cD} V_D^{\otimes n_D}$ as its Cartan
component. This yields morphisms of $G$-modules
$$
p_{D_1,D_2}^{\lambda} : V_{D_1} \otimes V_{D_2} \to 
\Sym(\bigoplus_{D \in \cD} V_D),
$$
obtained as the composition
$$
V_{D_1} \otimes V_{D_2} \to V_{D_1}V_{D_2} \to V(\lambda) \to
\bigotimes_{D\in \cD} V_D^{\otimes n_D} \to 
\Sym(\bigoplus_{D \in \cD} V_D).
$$
Note that each $p_{D_1,D_2}^{\lambda}$ is only defined up to
multiplication by a non-zero scalar, since the same holds for the
projection $V_{D_1} \, V_{D_2} \to V(\lambda)$.

\begin{proposition}\label{prop:rel}
With the above notation, there exist unique normalizations of the
morphisms $p_{D_1,D_2}^{\lambda}$ such that the ideal $I$ is generated
by the elements
\begin{equation}\label{eqn:rel}
v_1 \otimes v_2 - 
\sum_{\lambda \in \Lambda^+, \, \lambda \leq_X \omega_{D_1} + \omega_{D_2}} 
p_{D_1,D_2}^{\lambda}(v_1 \otimes v_2) \, \prod_{i=1}^r s_i^{n_i},
\end{equation}
where $v_1 \in V_{D_1}$, $v_2 \in V_{D_2}$, and 
$n_1,\ldots, n_r$ are the non-negative integers defined by
$\omega_{D_1} + \omega_{D_2} - \lambda = 
\sum_{i=1}^r n_i \, \gamma_i$.
\end{proposition}

\begin{proof}
It follows from Theorem \ref{thm:rees} that $I$ contains the elements
(\ref{eqn:rel}) for suitable normalizations of the
$p_{D_1D_2}^{\lambda}$. Let $J\subseteq I$ be the ideal of $S$
generated by all these elements. The natural map $S \to R(X)$ factors
through a surjective homomorphism of $G$-algebras  
$$
h : S/J \to R(X).
$$

We claim that the induced (surjective) homomorphism
$$
\bar{h}: S/(J,s_1,\ldots,s_r) \to R(X)/(s_1,\ldots,s_r)
$$
is an isomorphism. Indeed, the left-hand side is generated as a
$G$-algebra by simple modules $V_D$, $D \in \cD$, such that the
product of any two such modules is their Cartan product, or zero. By
a result of Kostant (see, e.g., \cite[4.1~Lemme]{Br85}), it follows
that the product of any two simple modules in this $G$-algebra is
again their Cartan product, or zero. Thus, the subalgebra 
$$
(S/(J,s_1,\ldots,s_r))^U = S^U/(J^U,s_1,\ldots,s_r)
$$
is generated by the images of the $U$-invariant lines in the $V_D$. On
the other hand, 
$$
(R(X)/(s_1,\ldots,s_r))^U = R(X)^U/(s_1,\ldots,s_r)
$$
is a polynomial ring in the images of the $s_D$, $D \in \cD$, by
Proposition \ref{prop:quot}. This implies our claim. 

From this claim and the freeness of the $k[s_1,\ldots,s_r]$-module
$R(X)$ (Proposition \ref{prop:quot} again), it follows that 
$$
\ker(h) = (s_1,\ldots,s_r) \ker(h).
$$ 
On the other hand, the algebra $S$ is graded by the weight group
$\Lambda$, where each $s_i$ has degree $\gamma_i$ and each $s_D$ has
degree $\omega_D$. Note that the degrees of non-zero homogeneous
elements of $S$ are linear combinations of simple roots with
non-negative coefficients, and that the homogeneous elements of degree
$0$ are just the constants. Hence $S$ admits a coarser positive
grading. Moreover, the generators of $J$ are all homogeneous, and the
ideal of $R(X)$ is homogeneous as well (the induced grading on $R(X)$
arises from the morphism 
$u : T \to \bbG_m^{\cD}$, $t \mapsto (\omega_D(t))$.) Thus, $\ker(h)$
is also graded. By the graded Nakayama lemma, it follows that 
$\ker(h) = 0$.
\end{proof}
 
\begin{proposition}\label{prop:gs}
For any spherical root $\gamma_i$, there exist colors $D_1,D_2$ and a
dominant weight $\lambda$ such that the $G$-module $V(\lambda)$
occurs in the product $V_{D_1} \, V_{D_2} \subset k[G/K]$, and 
$\omega_{D_1} + \omega_{D_2} - \lambda$ is a positive integral multiple
of $\gamma_i$.
\end{proposition}

\begin{proof}
Let $J$ be the ideal of the polynomial ring $k[s_1,\ldots,s_r]$
generated by those monomials $\prod_{i=1}^r s_i^{n_i}$ that appear in
some relation (\ref{eqn:rel}). Then  the quotient $R(X)/J R(X)$ is a
$G$-algebra generated by the images of the trivial $G$-modules
$s_1,\ldots,s_r$ and of the $G$-modules $V_D$, $D \in \cD$, and the
product of any two such modules is their Cartan product. As in the
proof of Proposition \ref{prop:rel}, it follows that the product of
any two simple modules in the $G$-algebra $R(X)/JR(X)$ is their Cartan
product, or zero. So, for any point $x$ of the zero set $Z(J)$ in
$\bbA^r$, the fiber $\tX_x$ is horospherical. Thus, $Z(J)$ is just the
origin, and hence $J$ contains positive powers of $s_1,\ldots,s_r$.
\end{proof}

\begin{example}\label{ex:group7}
The algebra $R(\bar{G}_{\ad})$ is generated by $s_1,\ldots,s_r$ and the
$G\times G$-submodules 
$\End(V(\omega_1)),\ldots,\End(V(\omega_r))$. For 
$1 \leq i \leq j \leq r$, we have a canonical isomorphism
$$
V(\omega_i) \otimes V(\omega_j) \cong \bigoplus_{\lambda \in \Lambda^+}
N_{ij}^{\lambda} \otimes V(\lambda),
$$
where 
$N_{ij}^{\lambda} := 
\Hom^G(V(\lambda), V(\omega_i) \otimes V(\omega_j))$.  
Thus, the tensor product 
$\End(V(\omega_i)) \otimes \End(V(\omega_j))$
is the direct sum of the 
$\End(N_{ij}^{\lambda}) \otimes \End(V(\lambda))$. Moreover, the
product $\End(V(\omega_i)) \, \End(V(\omega_j))$ in $k[G]$ is obtained
by projecting each $\End(N_{ij}^{\lambda})$ to $k$ via the trace map.
This yields the relations in $R(\bar{G}_{\ad})$:
$$
\End(V(\omega_i)) \, \End(V(\omega_j)) \cong 
\bigoplus_{\lambda}  \End(V(\lambda))\, \prod_{i=1}^r s_i^{n_i},
$$
where the sum runs over those dominant weights $\lambda$ such that
$V(\lambda)$ occurs in $V(\omega_i) \otimes V(\omega_j)$, and 
$\omega_i + \omega_j -\lambda = \sum_{i=1}^r n_i \, \alpha_i$.
Since the $G$-module $V(2\omega_i -\alpha_i)$ occurs in 
$V(\omega_i) \otimes V(\omega_i)$, the root monoid is generated by the
simple roots $\alpha_1,\ldots,\alpha_r$.

In particular, if $G = \SL_n$ then $V(\omega_i) = \bigwedge^i k^n$ for 
$i=1,\ldots, n-1$, and  
$$
V(\omega_i) \otimes V(\omega_j) \cong 
\bigoplus_{i',j', \; 0 \leq i' \leq i \leq j \leq j' \leq n} 
V(\omega_{i'}+\omega_{j'}),
$$
where we set $\omega_0 = \omega_n =0$. Thus, the relations have degree
at most $2$ in the $\End(V(\omega_i))$. This also holds for 
$G = \Sp_{2n}$, but not for any other simple group.
\end{example}

\section{The total coordinate ring of a spherical variety}
\label{sec:spherical}

\subsection{The equivariant divisor class group}
\label{subsec:class}

Let $X$ be a spherical variety under the connected reductive
group $G$. Replacing $G$ with a finite cover, we may assume that 
$G = \bG \times C$, where $C$ is a torus (the connected center of $G$)
and $\bG$ is a simply-connected semi-simple group (the derived
subgroup of $G$). Then $k[G]$ is a UFD, and $k[G]^* = k[C]^*$ consists
of the scalar multiples of characters of $C$ (identified with
characters of $G$). Moreover, replacing $C$ with a quotient torus, we
may assume that $C$ acts faithfully on $X$.

We denote by $X^0 \subseteq X$ the open $G$-orbit and by
$X_1, \ldots , X_n$ the irreducible components of codimension $1$ of
the boundary $X \setminus X^0$; they form the (possibly empty) set of
\emph{boundary divisors}. Let $X^{\leq 1}$ be the union of $X^0$ and
all $G$-orbits of codimension $1$, i.e., of the open $G$-orbits in the
boundary divisors. Clearly, $X^{\leq 1}$ is an open $G$-stable subset
of $X$, and its complement has codimension $\geq 2$. Moreover,
since $X$ is a normal $G$-variety, the open subset $X^{\leq 1}$ is
non-singular, and $k[X] = k[X^{\leq 1}]$. We shall assume for
simplicity that $k[X] = k$. This holds, e.g., if $X$ is complete. 

Recall that a \emph{divisorial sheaf} on $X$ is a coherent sheaf $\cF$
which is reflexive of rank $1$; equivalently, $\cF \cong \cO_X(D)$
for some Weil divisor $D$ on $X$. Then $\cF$ is canonically isomorphic
to $j_*(j^*\cF)$, where 
$$
j : X^{\leq 1} \to X
$$ 
denotes the inclusion; moreover, the sheaf $j^*\cF$ is
invertible. Conversely, any invertible sheaf $\cL$ on $X^{\leq 1}$
defines a divisorial sheaf $\cF := j_*\cL$, and $\cL$ may be
identified with $j^* \cF$. So we shall freely identify divisorial
sheaves on $X$ with invertible sheaves on $X^{\leq 1}$. In particular,
the tensor product of invertible sheaves defines the product of
divisorial sheaves (the double dual of their tensor product).

Together with \cite[Sec.~2.4]{KKLV}, it follows that any divisorial
sheaf has $G$-linearizations, and any two of them differ by a
unique character of $C$. The group of isomorphism classes of
$G$-linearized divisorial sheaves on $X$ may be identified with
$\Pic^G(X^{\leq 1})$; we denote this group by $\Cl^G(X)$ and call it
the \emph{equivariant divisor class group}. (This is the equivariant
Chow group $A_{\dim(X)-1}^G(X)$ of \cite{EG98}). Note the exact
sequence
\begin{equation}\label{eqn:clex}
0 \to \cX(C) \to \Cl^G(X) \to \Cl(X) \to 0,
\end{equation}
where $\Cl(X)$ denotes the usual class group. 

We shall obtain a presentation of the abelian group $\Cl^G(X)$. For
this, we generalize some of the notation of \ref{subsec:wonderful}: 
we define the open $B$-orbit $X^0_B$, the set $\cD$ of irreducible
components of $X^0 \setminus X^0_B$ (these are prime $B$-stable
divisors that we also regard as divisors in $X$), the weight group 
$\Lambda(X) \cong k(X)^{(B)}/k^*$, the valuations 
$v_{X_i}$, $v_D$ of the function field $k(X)$ associated with the
boundary divisors and the colors, and their restrictions
$v_i$, $\rho(v_D)$ to $\Hom(\Lambda(X),\bbR)$. We choose a
base point $x \in X^0_B$ and denote by $H$ its stabilizer in $G$.

The pull-back of each $D \in \cD$ under the projection 
$\pi: G \to G/H \cong X^0$ 
is a $B\times H$-stable divisor, which has a unique equation 
$f_D \in k[G]$ such that $f_D(1) = 1$ and $f_D$ is invariant
under left multiplication by $C$. Then $f_D$ is a 
$B\times H$-eigenvector; let $(\omega_D,\chi_D)$ be its weight. Also,
let $V_D \subset k[G]$ be the span of the left $G$-translates of $f_D$.
The $G$-module $V_D$ is isomorphic to $V(\omega_D)$ and consists
of right $H$-eigenvectors of weight $\chi_D$. It yields a
$G$-morphism $G/H \to \bbP(V_D)^*$ which extends to a $G$-morphism 
$$
\varphi : X^{\leq 1} \to \bbP(V_D)^*.
$$ 
Moreover, $D$ is the pull-back under $\varphi$ of the hyperplane of
$\bbP(V_D)^*$ orthogonal to $f_D \in V_D$. This yields a linearization
of the invertible sheaf 
$\cO_{X^{\leq 1}}(D) \cong \varphi^*\cO_{\bbP(V_D^*)}(1)$, 
and hence a \emph{canonical} linearization of the divisorial sheaf
$\cO_X(D)$. We denote by $[D]^G \in \Cl^G(X)$ the corresponding
class (notice, however, that $D$ is never $G$-stable). We also obtain
classes $[X_i]^G$ associated with boundary divisors, since the sheaves 
$\cO_X(X_i)$ are canonically $G$-linearized.

\begin{proposition}\label{prop:presentation}
(1) The abelian group $\Cl^G(X)$ is defined by generators
$$
[X_1]^G, \ldots, [X_n]^G, \quad [D]^G ~(D \in \cD), 
\quad \chi \in \cX(C)
$$
and relations 
\begin{equation}\label{eqn:clrell}
\sum_{i=1}^n v_i(\lambda) \, [X_i]^G +
\sum_{D \in \cD} \rho(v_D)(\lambda) \, [D]^G = \lambda \vert_C
\end{equation}
for all $\lambda \in \Lambda(X)\subseteq \cX(T)$.

\smallskip

\noindent
(2) Likewise, the abelian group $\Cl(X)$ is defined by generators
$[X_1]$, $\ldots$ , $[X_n]$, $[D]$ $(D\in \cD)$, and relations 
$$
\sum_{i=1}^n v_i(\lambda) \, [X_i] + 
\sum_{D \in \cD} \rho(v_D)(\lambda) \, [D] =0
$$ 
for all $\lambda \in \Lambda(X)$.
\end{proposition}

\begin{proof}
We begin with the easy proof of (2). Since $X^0_B$ has trivial
divisor class group, the group $\Cl(X)$ is generated by the classes of
boundary divisors and of colors. Moreover, all relations between these
$B$-stable divisors arise from divisors of rational functions on $X$
which are $B$-eigenvectors.

Together with (\ref{eqn:clex}), it follows that $\Cl^G(X)$ is
generated by the above classes and the characters of $C$. Moreover,
any $f \in k(X)^{(B)}$ may be regarded as a $C$-invariant rational
section of the trivial line bundle on $X$ with linearization twisted
by $-\lambda\vert_C$, where $\lambda$ denotes the weight of $f$. Thus,  
$\div(f) - \lambda\vert_C = 0$ in $\Cl^C(X) \cong \Cl^G(X)$. In other
words, the relations (\ref{eqn:clrell}) hold in $\Cl^G(X)$. 

Assume that 
$$
\sum_{i=1}^n x_i \, [X_i]^G + \sum_{D \in \cD} y_D \, [D]^G = \chi
$$
in $\Cl^G(X)$, where $x_i$, $y_D$ are integers and $\chi \in \cX(C)$.
Then 
$$
\sum_{i=1}^n x_i \, [X_i] + \sum_{D \in \cD} y_D \, [D] =0
$$ 
in $\Cl(X)$. Thus, there exists $\lambda \in \Lambda(X)$ such that
$v_i(\lambda) = x_i$ for all $i$, and $\rho(v_D)(\lambda) = y_D$ for
all $D$. Hence $\lambda\vert_C = \chi$ in $\Cl^G(X)$. By
(\ref{eqn:clex}), $\lambda\vert_C = \chi$ in $\cX(C)$, so that 
(\ref{eqn:clrell}) gives all the relations.
\end{proof}

As in \ref{subsec:picard}, the inclusion 
\begin{equation}\label{eqn:i}
i : G/H \to X
\end{equation}
defines the restriction map
$$
i^* : \Cl^G(X)  = \Pic^G(X^{\leq 1}) \to \Pic^G(G/H) \cong \cX(H),
$$
where the latter isomorphism is defined by $\chi \mapsto \cL(\chi)$. 
Moreover, the first assertion of Proposition \ref{prop:sequ} may be
easily generalized:

\begin{lemma}\label{lem:sequ-bis}
(1) The sequence
$$
\CD
0 @>>> \bbZ^n @>{u}>> \Cl^G(X) @>{i^*}>> \Pic^G(G/H) @>>> 0 \\ 
\endCD
$$
is exact, where $u(x_1,\ldots,x_n) = \sum_{i=1}^n x_i \, [X_i]^G$. 

\smallskip

\noindent
(2) The canonical section $s_D \in \Gamma(X,\cO_X(D))$ is a
$B$-eigenvector of weight $\omega_D$; the canonical section 
$s_i \in \Gamma(X,\cO_X(X_i))$ is $G$-invariant.
\end{lemma}

\begin{example}
In the case where $X$ is a toric variety (i.e., $G = T = C$ is a
torus), the group $\Cl^T(X)$ is freely generated by the classes of
boundary divisors, as follows, e.g., from Lemma \ref{lem:sequ-bis}.  
In particular, this group is torsion-free, while the torsion subgroup
of $\Cl(X)$ may be any prescribed finite abelian group. 

For an arbitrary spherical $G$-variety $X$, the group $\Cl^G(X)$ may
have non-zero torsion. Specifically, consider the standard embedding
$X$ of the group $\SO_{2n}$, i.e., the normalization of the wonderful
compactification of the semi-simple adjoint group $\PSO_{2n}$ in the
function field of $\SO_{2n}$. Then 
$G = \bG = \Spin_{2n} \times \Spin_{2n}$, so that $\Cl^G(X) = \Cl(X)$. 
Using Proposition \ref{prop:presentation}, one checks that this group
is the direct product of its subgroup $\Pic(X)$ (freely generated by
the classes of the $n$ colors) and a group of order $2$.
\end{example}

Some of the results of \ref{subsec:cone} also extend to spherical
varieties. Specifically, let $\Eff^G(X)\subseteq \Cl^G(X)$ be the
submonoid consisting of classes of sheaves having non-zero global
sections. Then $\Eff^G(X)$ is the pull-back under the map
$\Cl^G(X) \to \Cl(X)$ of the effective monoid
$\Eff(X) \subseteq \Cl(X)$, and the latter is generated by the classes
of boundary divisors and of colors. Moreover, the associated cone
$\Eff(X)_{\bbR} \subset \Cl(X)_{\bbR}$ is a rational polyhedral cone,
which contains no line since $k[X] = k$. 

The extremal rays of $\Eff(X)_{\bbR}$ which contain no color are still
the rays generated by the fixed boundary divisors, as in Lemma
\ref{lem:boundary}. Moreover, the extremal rays containing no boundary
divisor correspond to $G$-morphisms $X^{\leq 1} \to G/P$
where $P \supseteq H$ is a maximal parabolic subgroup of $G$, as in
Lemma \ref{lem:color}. But the statement of Lemma \ref{lem:inv}
does not hold in general, since two boundary divisors may well have
the same class (e.g., if $G = \bbG_m$ and $X = \bbP^1$). Likewise,
Lemma \ref{lem:mixed} admits no straightforward generalization, as
shown by the example of the projective cone over $\bbP^1$ embedded via
the sections of $\cO(n)$.

We shall determine the extremal rays of $\Eff(X)_{\bbR}$ when $X$ is
the \emph{standard embedding} of a \emph{sober} spherical homogeneous
space $G/H$, i.e., the group $N_G(H)/H$ is finite, and $X$ is the
unique complete toroidal embedding with only one closed orbit. 
Note that $H$ contains $C$, so that we may assume $G$ to be
semi-simple. The variety $X$ is generally singular, but with quotient
singularities only; it admits a finite surjective $G$-morphism 
$$
f: X \to \bX,
$$ 
where $\bX$ is a wonderful $G$-variety, and the preimage of any color
of $\bX$ is a unique color (see \cite[Cor.~7.6]{Kn96}; in fact, there
is a canonical choice of $\bX$, see \cite{Lu01} or Subsection
\ref{subsec:structure}). Moreover, the Picard group of $X$ is still
freely generated by the classes of colors, see \cite{Br89b}. Thus,
$f^* : \Pic(\bX) \to \Pic(X)$ yields isomorphisms 
$$
\Pic(\bX)_{\bbR} \cong \Pic(X)_{\bbR} = \Cl(X)_{\bbR} = \Cl^G(X)_{\bbR},
\quad \Eff(\bX)_{\bbR} \cong \Eff(X)_{\bbR}.
$$
Recall from Theorem \ref{thm:moved} that the extremal rays of
$\Eff(\bX)_{\bbR}$ are generated by the classes of fixed boundary
divisors, and of colors associated with surjective morphisms 
$\bX \to \bX'$ with connected fibers, where $\bX'$ is a wonderful
$G$-variety of rank $\leq 1$ having a unique color. Moreover, $f$
induces a bijection between the boundary divisors 
$\bX_1, \ldots, \bX_r$ and  $X_1, \ldots, X_r$, 
which restricts to a bijection between the subsets of fixed divisors
(since $f$ is finite). Likewise, the morphisms $\bX \to \bX'$ as above
are in bijection with morphisms $X \to X'$ satisfying the same
assumptions. This yields: 

\begin{proposition}\label{prop:sober}
Let $X$ be the standard embedding of a sober spherical homogeneous
space. Then the extremal rays of $\Eff(X)_{\bbR}$ consist of the rays
generated by fixed boundary divisors, and of those generated by
pull-backs of colors under $G$-morphisms $\varphi: X  \to X'$, where
$X'$ is a wonderful $G$-variety of rank $1$ having a unique color (in
particular, $X'$ is homogeneous). Both types of rays are disjoint.
\end{proposition}

Also, note that a boundary divisor is either fixed or semi-ample
(i.e., some positive multiple is base-point-free). However, a
semi-ample boundary divisor $X_i$ may well satisfy 
$\Gamma(X,\cO_X(X_i)) = k s_i$ and hence be stable under
$\Aut^0(X)$. This happens, e.g., for the two boundary divisors of the 
standard embedding of the group $\SO_4$.

\subsection{The total coordinate ring}
\label{subsec:total}

Our first aim is to define this ring $R(X)$, where $X$ is a spherical
variety under the connected reductive group $G$. The idea is to choose
$G$-linearized divisorial sheaves $\cF(\delta)$ in all the classes
$\delta \in \Cl^G(X)$, and to endow the sum 
$$
\bigoplus_{\delta \in \Cl^G(X)} \Gamma(X, \cF(\delta))
$$
with a multiplication law. For this, we need morphisms
$$
\cF(\delta) \otimes \cF(\delta') \to \cF(\delta + \delta')
$$
so that the resulting law is commutative and associative. The
existence of such representatives and morphisms is guaranteed by the
following: 

\begin{lemma}\label{lem:sections}
Every $\delta \in \Cl^G(X)$ is the class of a unique $G$-linearized
divisorial subsheaf 
$$
\cF(\delta) \subseteq i_* \cL(\chi),
$$
where $i$ is defined by (\ref{eqn:i}), and 
$\chi = i^*(\delta) \in \cX(H)$. In particular,
$\Gamma(X,\cF(\delta))$ is a $G$-submodule of $k[G]^{(H)}_{\chi}$
which depends only on $\delta$. 

Moreover, for any $\delta' \in \Cl^G(X)$ with associated character
$\chi' \in \cX(H)$, the canonical homomorphism 
$i_* \cL(\chi) \otimes i_* \cL(\chi') \to i_* \cL(\chi + \chi')$
maps $\cF(\delta) \otimes \cF(\delta')$ to $\cF(\delta + \delta')$.
\end{lemma}

\begin{proof}
Let $\cF$ be a representative of $\delta$. Then the natural map 
$\cF \to i_*(i^*\cF)$ identifies $\cF$ with a subsheaf of 
$i_* \cL(\chi)$. To show that this subsheaf only depends 
on $\delta$, we may assume that $X = X^{\leq 1}$. Then $X$ is covered
by $G$-stable open subsets consisting of two $G$-orbits: the open
orbit and an orbit of codimension $1$. So we may further assume that
$X$ is one of these ``elementary embeddings'' of $G/H$.

Let $P \supseteq B$ be the stabilizer in $G$ of the open $B$-orbit
$X^0_B$, and let $X_B \subseteq X$ be the complement of the union of
the colors. Then $X_B$ is an affine $P$-stable open subset of $X$,
isomorphic to $P\times^L \bbA^1$ where $L$ is a Levi subgroup of $P$,
acting on $\bbA^1$ via a character $\lambda$; this induces a
$P$-equivariant isomorphism $X^0_B \cong P \times^L \bbG_m$ (see,
e.g., \cite[Sec.~1.1]{Br89b}). Since $X$ is covered by the $G$-translates
of $X_B$, it suffices to check that the image of $\Gamma(X_B,\cF)$ in
$\Gamma(X^0_B,\cL(\chi))$ depends only on $\delta$. 

Using the correspondence between $P$-linearized sheaves on
$P\times^L \bbA^1$ and $L$-linearized sheaves on $\bbA^1$, we may
assume that $X = \bbA^1$ where $G = L$ acts via $\lambda \in \cX(C)$. 
Then $\cF$ is isomorphic to $\cO_X$ with linearization twisted
by $\chi \in \cX(C)$, and $\delta$ may be identified with $\chi$. 
Moreover, $\Gamma(X^0_B,\cL(\chi))$ is the subspace of $k[C]$ spanned
by the characters $\chi + n \lambda$ where $n \in \bbZ$, while the
image of $\Gamma(X,\cF)$ is spanned by the $\chi + n \lambda$ where
$n$ is non-negative. This completes the proof of the first assertion;
the second one may be checked similarly.
\end{proof}

As a straightforward consequence, we obtain the following: 

\begin{proposition}\label{prop:alg}
With the notation of Lemma \ref{lem:sections}, the sum 
$$
\bigoplus_{\delta \in \Cl^G(X)} \cF(\delta)
$$ 
has a natural structure of a sheaf of $\cO_X$-algebras, graded by
$\Cl^G(X)$, and compatibly $G$-linearized. In particular,
$$
R(X) := \bigoplus_{\delta \in \Cl^G(X)} \Gamma(X,\cF(\delta))
$$
has a natural structure of a $G$-algebra, graded by the abelian group
$\Cl^G(X)$, and depending only on $X^{\leq 1}$. 
\end{proposition}

Next, recall that $\cX(C)$ is a subgroup of $\Cl^G(X)$, and the
corresponding $G$-linearized sheaves are just copies of the structure
sheaf $\cO_X$ with linearization twisted by characters. Since
$\Gamma(X,\cO_X) = k$, the subalgebra
$$
\bigoplus_{\delta \in \cX(C)} \Gamma(X,\cF(\delta))\subseteq R(X)
$$
is isomorphic to the coordinate ring $k[C]$. Together with
(\ref{eqn:clex}), it follows that $R(X)^* \cong k[C]^* = k^* \cX(C)$.

By the same arguments as in Proposition \ref{prop:quot}, we obtain:

\begin{proposition}\label{prop:quotbis}
The invariant algebra $R(X)^U$ is freely generated as a 
$k[C]$-algebra by $s_1,\ldots,s_n$ and the $s_D$, $D \in \cD$.  

In particular, the subalgebra $R(X)^{\bG}$ (resp.~$R(X)^G$) is
freely generated by $s_1,\ldots,s_n$ as a $k[C]$-algebra
(resp.~as a $k$-algebra).
\end{proposition}

Together with \cite{Po87}, it follows that the $k$-algebra $R(X)$ is
finitely generated and normal. Thus, it corresponds to a normal affine
variety $\tX$ equipped with an action of the group 
$$
\tG := G \times \Gamma_X,
$$
where we put 
$$
\Gamma_X := \Hom(\Cl^G(X),\bG_m),
$$ 
the diagonalizable group with character group $\Cl^G(X)$. The sheaf of
$\cO_X$-algebras $\bigoplus_{\delta \in \Cl^G(X)} \cF(\delta)$ yields
a $\tG$-variety equipped with an affine morphism $p:\hX \to X$ which
is $G$-equivariant and $\Gamma_X$-invariant. Moreover, $p$ is a good 
quotient for the action of $\Gamma_X$ (that is, the natural map
$\cO_X \to p_*(\cO_{\hat X})^{\Gamma_X}$ is an isomorphism), and the
restriction of $p$ to the open subset $X^{\leq 1}$ (or to the regular
locus of $X$) is a principal bundle for that group.

As in \cite{BH03} or in \ref{subsec:basic}, one checks that the
natural map $\iota : \hX \to \tX$ is an open immersion, the
complement of its image has codimension $\geq 2$, and $R(X)$ is a
UFD. Moreover, $\hX$ contains an open orbit of a Borel subgroup
of the connected reductive group $\tG^0$, so that the $\tG^0$-variety
$\tX$ is spherical.

One may show, as in Proposition \ref{prop:quot}, that the quotient map
$\tX \to \tX/\!/ G$ is flat, and its fibers are normal varieties. This
also follows from a structure theorem for $R(X)$ that we shall obtain
in the next subsection.

\begin{example}
In the toric case, Proposition \ref{prop:quotbis} yields an algebra
isomorphism $R(X) \cong k[T][s_1,\ldots,s_n]$; equivalently,
$\tX \cong T \times \bbA^n$. Moreover, $\tG \cong T \times \bbG_m^n$,
where $T$ (resp.~$\bbG_m^n$) acts on $\tX$ via its natural action on
the first (resp.~second) factor. One may check that the complement 
$\tX \setminus \hX$ is the union of the subspaces 
$(x_i=0),~i\in I$, where the $x_i$ are the coordinates on $\bbA^n$,
and $I \subseteq \{1,\ldots,n\}$ is such that $\bigcap_{i \in I} X_i$
is empty.
\end{example}

\subsection{Structure}
\label{subsec:structure}

We begin by recalling how to associate a wonderful variety to
any spherical variety, after \cite{Lu01}. The notation is as in 
\ref{subsec:class}; in particular, $X$ denotes a spherical variety
under $G = \bG \times C$, with open orbit $G/H$.

The normalizer $N_G(H)$ acts on $G/H$ by equivariant automorphisms,
and permutes the colors $D \in \cD$. The subgroup of $N_G(H)$ which
stabilizes all colors is called the \emph{spherical closure} of
$H$ in \cite{Lu01}. Clearly, this subgroup has finite index in
$N_G(H)$, and contains both groups $H$ and $C$; thus, it may be
written as $\bH \times C$, where $\bH$ is a closed subgroup of the
derived subgroup $\bG$. The homogeneous space $\bG/\bH$ is spherical
and sober; its set of colors may be identified with $\cD$. 

We record the following observation, implicit in \cite[Sec.~6]{Lu01}.

\begin{lemma}\label{lem:kernel}
Let $K$ be the intersection of the kernels of all characters of
$H$. Then $K$ is also the intersection of the kernels of all
characters of $\bH$.
\end{lemma}

\begin{proof}
Denote the latter group by $\bK$, we may identify this group with the
intersection of all characters of $\bH \times C$. Then 
$\bK \subseteq H$ and $H/\bK$ is diagonalizable, since the quotient
$(\bH \times C)/H \subseteq N_G(H)/H$ is diagonalizable. Thus, 
$K \subseteq \bK \subseteq H$.

The group $\cX(H)$ is generated by the restrictions of characters of
$C$ and the weights of the $f_D$, $D \in \cD$ (see, e.g., the proof of
Lemma 6.3.1 in \cite{Lu01}). Thus, each character of $H$ extends to a
character of $\bH \times C$. It follows that $\bK \subseteq K$.
\end{proof}

By \cite[Cor.~7.6]{Kn96}, the homogeneous space $\bG/\bH$ admits a
unique wonderful embedding $\bX$. The natural map
$$
G/H \to G/(\bH \times C) = \bG/\bH
$$
extends uniquely to a $G$-morphism
$$
\varphi : X^{\leq 1} \to \bX
$$
which is $C$-invariant; moreover, $\varphi$ extends to the whole
variety $X$ if and only if this variety is toroidal. 

Recall that each invertible sheaf $\cL$ on $\bX$ carries a canonical
$\bG$-linearization; thus, its pull-back $\varphi^*\cL$ is
$G$-linearized, where $C$ acts trivially. This defines a homomorphism
\begin{equation}\label{eqn:pb}
\varphi^* : \Pic(\bX) = \Pic^{\bG}(\bX) \to 
\Pic^G(X^{\leq 1}) = \Cl^G(X)
\end{equation}
and morphisms of $G$-modules
$$
\varphi^*: \Gamma(\bX,\cL) \to \Gamma(X^{\leq 1},\varphi^*\cL)
$$
(via the natural map 
$\cL \to \varphi_*(\varphi^*\cL) 
= \cL \otimes \varphi_* \cO_{X^{\leq 1}}$). Note that 
\begin{equation}\label{eqn:pbc}
\varphi^*([\bD]) = [D]^G
\end{equation} 
for each $D \in \cD$, where $\bD$ denotes the corresponding color of
$\bX$.
 
In turn, we obtain a homomorphism of graded $G$-algebras 
$$
\varphi^*: R(\bX) \to R(X),
$$
where the corresponding homomorphism between their grading groups is
given by (\ref{eqn:pb}). 

\begin{theorem}\label{thm:cart}
(1) The restriction
$$
\varphi^{*\bG}: R(\bX)^{\bG} = k[\bs_1,\ldots,\bs_r] \to 
R(X)^{\bG} = k[C][s_1,\ldots,s_n]
$$
is given by 
$$
\varphi^{*\bG}(\bs_i) = \prod_{j=1}^n s_j^{-v_j(\gamma_i)}
$$
for $i=1,\ldots,r$.

\smallskip

\noindent
(2) The map
$$
\psi: R(\bX) \otimes_{R(\bX)^{\bG}} R(X)^{\bG} \to R(X),
\quad
\bs \otimes s \mapsto \varphi^*(\bs) s
$$
is an isomorphism.

\smallskip

\noindent
(3) The multiplication map induces an isomorphism
$$
k[C] \otimes R(X)^C \cong R(X).
$$
Moreover, the $G$-module $R(X)^C$ is generated by the monomials in
$s_1,\ldots,s_r$ and the $s_D$, $D \in \cD$. 
\end{theorem} 

\begin{proof}
(1) Recall that 
$\div_{\bX}({\bf f}_i) = \bX_i - \sum_{D \in \cD} 
\langle D,\gamma_i \rangle \, \bD$
with obvious notation. Thus, 
$\div_X(\varphi^*{\bf f}_i) = \varphi^*(\bX_i) - \sum_{D \in \cD} 
\langle D,\gamma_i \rangle \, D$.
It follows that 
$\varphi^*[\bX_i] = \sum_{j=1}^n v_j({\bf f}_i) [X_j]^G$,
which implies our assertion.

(2) It suffices to show that the restriction
$$
\psi^U: R(\bX)^U \otimes_{R(\bX)^{\bG}} R(X)^{\bG} \to R(X)^U
$$
is an isomorphism. The left-hand side is generated as a $k[C]$-algebra 
by the $\bs_{\bD} \otimes 1$, $D\in \cD$, and the $1 \otimes s_i$,
$i=1,\ldots,n$. Moreover, $\psi(\bs_{\bD} \otimes 1) = s_D$ by
(\ref{eqn:pbc}), and $\psi(1 \otimes s_i) =s_i$. So the assertion
follows from Proposition \ref{prop:quotbis}.

(3) follows from (2) together with the isomorphism 
$k[C] \otimes R(X)^G \cong R(X)^{\bG}$
given by the multiplication map. 
\end{proof}

In geometric terms, the quotient map $\tX \to \tX/\!/\bG$ is the
pull-back of the corresponding map 
$$
\bq : \tbX \to \tbX/\!/\bG \cong \bbA^r
$$ 
under the natural morphism 
$$
\tvarphi/\!/\bG : \tX/\!/\bG \to \tbX/\!/\bG.
$$ 
Moreover, $\tX/\!/\bG \cong C \times \tX/\!/G$ and $\tvarphi/\!/\bG$
factors through the second projection $\tX/\!/\bG \to \tX/\!/G$. It
follows that the quotient morphism
$$
q: \tX \to \tX/\!/G \cong \bbA^n
$$ 
is flat, and its fibers are exactly the products of $C$ with the
fibers of $\bq$.

Also, $\tX$ is a trivial principal $C$-bundle over 
$\tX/C := \Spec R(X)^C$. We may regard $R(X)^C$ as the total
coordinate ring of $X$ with respect to the canonical linearizations
defined in \ref{subsec:class}, so that the principal bundle 
$\tX \to \tX/C$ corresponds to the variation of linearizations.

\subsection{The automorphism group of a complete toroidal variety} 
\label{subsec:automorphisms-bis}

We shall extend part of the results of \ref{subsec:automorphisms} to
the setting of a complete toroidal variety $X$. We first relate the
connected automorphism group of $X$ (a linear algebraic group) to that
of a standard embedding, as in \ref{subsec:class}. Let 
$\varphi: X \to \bX$ be the morphism defined in \ref{subsec:structure}
and consider its Stein factorization
$$
\CD
X @>{\psi}>> X' @>{f}>> \bX.
\endCD
$$ 
Then $X' = \Spec_{\bX}(\varphi_*\cO_X)$ is a complete spherical
$G$-variety and $\psi$, $f$ are equivariant. Let $H'$ be the
stabilizer of the base point $\psi(x)$ of $X'$. Then
$$
C \cdot H \subseteq H' \subseteq C \times \bH \subseteq N_G(H)
$$
and $H'/H = (N_G(H)/H)^0$ is a torus containing $C$, while 
$N_G(H')/H' = N_G(H)/H'$ is a finite abelian group. The action of
$H'/H$ on $G/H$ extends to an action on $X$, so that 
$$
H'/H = \Aut^0_G(X).
$$ 
Since the morphism $f$ is finite, $X'$ is the standard embedding of
$G/H'$. Moreover, $\psi^{-1}(G/H') \cong G \times^{H'} F$, where $F$ is
a complete toric variety under $H'/H$; note that $F$ is the general
fiber of the morphism $\psi$.

The group $\Aut^0(X)$ acts on $X'$ so that $\psi$ is equivariant (this
follows, e.g., by considering the pull-back $L$ of an ample invertible
sheaf on $X'$, and linearizing a positive power of $L$ with respect to
$\Aut^0(X)$). This yields a homomorphism of algebraic groups 
$$
\psi_*: \Aut^0(X) \to \Aut^0(X')
$$
with kernel a closed subgroup of $\Aut(F)$ containing $H'/H$. The
structure of $\Aut(F)$ is well-known; in particular, $H'/H$ is a
maximal torus of this group, and hence of the kernel of $\psi_*$.
(However, there may exist no homomorphism
$\varphi_* : \Aut^0(X) \to \Aut^0(\bX)$, e.g., when $G = \Spin_{n+1}$
where $n\geq 3$, and $X \subset \bbP^{n+1}$ is a non-singular quadric
hypersurface. Then $X' = X$ and $\bX = \bbP^n$, so that
$\Aut^0(X) \cong \PSO_{n+2}$ while $\Aut^0(\bX) \cong \PSL_{n+1}$.)

We may now state our generalization of Theorem \ref{thm:aut}(2): 

\begin{theorem}\label{thm:aut-bis}
(1) The closed connected subgroup 
$\Aut^0(X, \partial X) \subseteq \Aut^0(X)$
is independent of the complete toroidal embedding $X$ of $G/H$.
The Lie algebra of this group is $\Gamma(X,\cS_X)$, where $\cS_X$
denotes the direct image of the logarithmic tangent sheaf of 
$X^{\leq 1}$. We have a split exact sequence of Lie algebras 
\begin{equation}\label{eqn:exact}
\CD
0 @>>> \fh'/\fh @>>> \Gamma(X,\cS_X) @>{\psi_*}>> 
\Gamma(X',\cS_{X'}) \to 0
\endCD
\end{equation}
and an isomorphism of Lie algebras 
\begin{equation}\label{eqn:isom}
\Gamma(X',\cS_{X'}) \cong \Gamma(\bX, \cS_{\bX}).
\end{equation}
In particular, the group $\Aut^0(X,\partial X)$ is reductive, with
connected center $H'/H$.

\smallskip

\noindent
(2) Let $G' \subseteq \Aut^0(X)$ be a closed connected subgroup
containing the image of $G$, and centralizing $H'/H$. Then $G'$ is 
reductive, its connected center is contained in $H'/H$, and the
$G'$-variety $X$ is toroidal, with the same colors as those of the
$G$-variety $X$. 
\end{theorem}

\begin{proof}
(1) We may choose a desingularization 
$$
\pi: Y \to X
$$ 
which is equivariant with respect to $\Aut^0(X)$. Then $Y$ is a
complete toroidal embedding of $G/H$, and we have an injective
homomorphism $\Aut^0(Y) \to \Aut^0(X)$ which restricts to 
$$
\pi_*: \Aut^0(Y,\partial Y) \to \Aut^0(X,\partial X).
$$ 
On the other hand, $\Aut^0(X,\partial X)$ lifts to a subgroup of
$\Aut^0(Y)$, which preserves the open $G$-orbit and hence 
$\partial Y$. Thus, $\pi_*$ is an isomorphism.

Let $\cS_Y$ be the sheaf of vector fields on $Y$ which preserve the
boundary. Then the sheaf $\pi_*(\cS_Y)$ is equipped with a map to
$\cS_X$. Using the local structure of toroidal varieties as in the
proof of \cite[Prop.~2.7.1]{BB96}, one checks that this map is an
isomorphism of locally free sheaves (indeed, one reduces to the toric 
case, where $\pi_*(\cS_Y)= \cS_X = \cO_X \otimes \fg$). Moreover, the
Lie algebra of $\Aut^0(Y,\partial Y)$ is $\Gamma(Y,\cS_Y)$ and does
not depend on the complete toroidal non-singular embedding $Y$ of
$G/H$; see, e.g., \cite[Sec.~2]{BB96}. Thus, the Lie algebra
$\Gamma(X,\cS_X) \cong \Gamma(Y,\cS_Y)$ 
depends only on $G/H$, so that the same holds for 
$\Aut^0(X,\partial X)$ (regarded as a group of automorphisms of
$G/H$). In fact, the $G$-variety $X$ is pseudo-free in the sense of
\cite[Sec.~2]{Kn94}, and the sheaf $\cS_X$ is the same as the one
defined there.
 
Likewise, the local structure of toroidal varieties yields an exact
sequence of locally free sheaves
$$
0 \to \cO_X \otimes (\fh'/\fh) \to \cS_X \to \psi^*(\cS_{X'}) \to 0.
$$
Taking global sections yields the exact sequence (\ref{eqn:exact}),
since $\psi_*(\cO_X) = \cO_{X'}$ and $H^1(X,\cO_X)=0$. On the other
hand, the isomorphism (\ref{eqn:isom}) follows from
\cite[Cor.~6.2]{Kn94}. Together with Theorem \ref{thm:aut}, this
implies that the Lie algebra $\Gamma(X',\cS_{X'})$ is semi-simple. In
turn, this implies that $\Aut^0(X,\partial X)$ is reductive with
connected center $H'/H$ (this may also be proved directly, along the
lines of Lemma \ref{lem:ss}).

(2) The kernel of the restriction of $\psi_*$ to $G'$ is a subgroup of
$\Aut(F)$ centralized by $H'/H$, i.e., a subgroup of $H'/H$. Thus, it
suffices to show that $\psi_*(G') \subseteq \Aut^0(X')$ is
semi-simple, and the $\psi_*(G')$-variety $X'$ is toroidal. In other
words, we may assume that $X$ is a standard embedding. Recall that 
$f: X \to \bX$ induces a bijection between the boundary divisors 
of $X$ and those of $\bX$, which preserves the fixed divisors

We claim that a boundary divisor $X_i$ is either $G'$-stable or
contains no $G'$-orbit. Indeed, if $X_i$ is moved by $G'$, then it is
not fixed, and hence semi-ample. This yields a surjective morphism
$\varphi: X \to X'$ with connected fibers, where $X'$ is a wonderful
variety of rank $1$ with ample boundary divisor $X'_i$, and 
$X'_i = \varphi^{-1}(X_i)$ (as sets). The group $G'$ acts on $X'$ so
that $\varphi$ is equivariant, and moves $X'_i$. Hence $G'$ acts
transitively on $X'$, which implies our claim.

Using that claim, we may argue as in the proof of Theorem
\ref{thm:aut}: Let $I$ be the set of indices of $G'$-stable boundary
divisors (so that $I$ contains the indices of fixed divisors) and let
$X_I := \bigcap_{i \in I} X_i$. Then $X_I$ is $G'$-stable, and any proper
$G$-stable subvariety contains no $G'$-orbit. Thus, $X_I$ is a unique
orbit of $G'$, and hence of a Levi subgroup $G''\subseteq G'$
containing $G$. By the claim again, the $G''$-variety $X$ is toroidal
with boundary divisors $X_i$, $i \in I$. Moreover, 
$G'' \subseteq G' \subseteq \Aut^0(X,\partial_I X)$. 
It follows that $G'$ is semi-simple, by the first step together with
Lemma \ref{lem:ss}. 

The assertion on colors is also checked by adapting the proof of
Theorem \ref{thm:aut}.
\end{proof}

\begin{corollary}
Let $X$ be a complete spherical $G$-variety and denote by
$\Aut_{\rm cr}(X)$ the largest closed connected subgroup of $\Aut^0(X)$
which stabilizes all the $G$-orbits. 

Then $\Aut_{\rm cr}(X)$ is reductive and depends only on the open
$G$-orbit in $X$. 

If $X$ is toroidal, then $\Aut_{\rm cr}(X) = \Aut^0(X,\partial X)$. 
\end{corollary}

\begin{proof}
Consider the rational $G$-equivariant map $f : X  - \to \bX$
as in 4.3. By \cite{RY02}, there exists a sequence of blowups
$\pi : Y \to X$ with nonsingular $G$-stable centers, such that the
composition $f \circ \pi : Y \to \bX$ is regular. Thus, $Y$ is a
complete toroidal $G$-variety, and the action of $\Aut_{\rm cr}(X)$ on
$X$ lifts to an action on $Y$. Clearly, the latter action preserves
the boundary $\partial Y$, and hence we may regard $\Aut_{\rm cr}(X)$
as a subgroup of $\Aut^0(Y,\partial Y)$. 

On the other hand, $\Aut^0(Y,\partial Y)$ preserves all the $G$-orbits
in $Y$, since their closures are exactly the partial intersections of
boundary divisors. So, via the homomorphism $\Aut^0(Y) \to \Aut^0(X)$
induced by $\pi$, the subgroup $\Aut^0(Y,\partial Y)$ is mapped onto 
$\Aut_{\rm cr}(X)$. We have shown that 
$\Aut_{\rm cr}(X) = \Aut^0(Y,\partial Y)$. By Theorem \ref{thm:aut-bis},
this group is reductive and depends only on the open $G$-orbit.
\end{proof}

Via this corollary, we may associate with any spherical homogeneous
space $G/H$ a connected reductive group $\Aut_{\rm cr}(G/H)$ that we
call the \emph{group of completely regular automorphisms}, by analogy
with the completely regular differential operators of \cite{Kn94}. 
 
It would be interesting to determine the spherical homogeneous spaces
$G/H$ whose group of completely regular automorphisms strictly
contains the image of $G$. By Theorem \ref{thm:aut-bis}, this
reduces to determining the wonderful varieties $X$ such that
$\Aut^0(X,\partial X)$ strictly contains the image of $G$. Those of
rank $0$ are exactly the flag varieties $G/P$ which are homogeneous
under a larger connected semi-simple group; their classification is
well-known, see \cite{De77}. Examples of rank $1$ include the quadrics
of dimension $6$ (where $G = \G_2$ and $\Aut^0(X,\partial X) = \SO_7$)
and of dimension $7$ (where $G = \Spin_7$ and 
$\Aut^0(X,\partial X) = \PSO_8$).

\end{document}